\documentstyle[12pt]{article}

\begin{document}
\title{Functions of several Cayley-Dickson variables and
manifolds over them.}
\author{S.V. L\"udkovsky}
\date{31 May 2004}
\maketitle
\begin{abstract}
Functions of several octonion variables are investigated and
integral representation theorems for them are proved. With the help
of them solutions of the ${\tilde {\partial }}$-equations are
studied. More generally functions of several Cayley-Dickson
variables are considered. Integral formulas of the
Martinelli-Bochner, Leray, Koppelman type used in complex analysis
here are proved in the new generalized form for functions of
Cayley-Dickson variables instead of complex. Moreover, analogs of
Stein manifolds over Cayley-Dickson graded algebras are defined and
investigated.
\end{abstract}
\section{Introduction}
\par In previous papers functions of one quaternion and
Cayley-Dickson variables were investigated \cite{luoyst,luoystoc}.
In them superdifferentiability of functions was studied and the
theory of holomorphic functions was investigated. It was done with
the help of line integration introduced and studied there. This line
integration restricted on complex functions gives ordinary Cauchy
integral, but ordinary Cauchy integral can not be in the usual
manner extended on continuous functions of Cayley-Dickson numbers.
\par This line integral is additive by rectifiable paths and
continuous functions in open domains in Cayley-Dickson algebras
${\cal A}_p$, it is also $\bf R$-homogeneous and left and right
linear over quaternions, but generally nonlinear relative to
quaternions or octonions or Cayley-Dickson algebras of higher order.
Over Cayley-Dickson algebras of higher order than quaternions left
or right linearity of operators certainly undermines relations
between ordered (associated) products of generators of
Cayley-Dickson numbers, for example, $(ij)l=-i(jl)$ for generators
of the octonion algebra $\bf O$, where $\{ 1,i,j,k \}$ are standard
generators of the quaternion algebra $\bf H$, $l$ is the generator
of doubling procedure of the construction of $\bf O$ from $\bf H$
\cite{baez,kurosh}. The line integral over Cayley-Dickson algebras
for $z$-superdifferentiable (that is, ${\cal A}_p$-holomorphic)
functions in corresponding domains $U$ depends only on specific
homotopy classes of rectifiable paths with given initial and final
points and also satisfies at least locally $(\partial [\int_{\gamma
_z} f(\xi )d\xi ]/\partial z).1 = f(z)$, for example, in a ball $B$
in $U$, where $\gamma _z: [0,1]\to U$, $\gamma (0)$ and $\gamma _z
(1)=z\in B$.
\par The Cayley-Dickson algebra is $\bf Z_2$-graded, that is, superalgebra.
In the theory of superalgebras it was traditionally used the notion
of superdifferentiabilty (left or right superlinear) \cite{khren}.
It causes strong restrictions on the types of admissible functions.
For example, over Grassman algebras it produces functions only
linear in odd arguments \cite{berez,dewitt}. In general it leads to
conditions analogous to Cauchy-Riemann.
\par Cayley-Dickson algebras such as quaternion and octonion algebras
have found applications in quantum mechanics and noncommutative
geometry \cite{connes,emch,lawmich,moore,oystaey}. The latter is
especially valuable in conjunction with operator algebras, which
permits to consider quantization. On the other hand, Cayley-Dickson
algebras are not central over the field of complex numbers $\bf C$.
Moreover, the octonion algebra and Cayley-Dickson algebras of higner
order can not be written as matrices with entries in the field of
real or complex numbers, though their centre is $\bf R$.
\par It is necessary to note, that apart from the real or complex
case derivatives of superdifferentiable functions of Cayley-Dickson
numbers are operators even in the case of one variable. Therefore, a
superderivative of the line integral by a final point of the path is
an operator. To work with rings of superdifferentiable functions it
was introduced the condition of ${\cal A}_p$-additivity instead of
right or left superlinearity of a superdifferential
\cite{luoyst,luoystoc}. It is natural, since if to start from the
family $\cal F$ of all right superlinearly superdifferentiable
functions $f: U\to {\cal A}_p$, $p\ge 2$, then the using of Leibnitz
rule for finite ordered (associated) products $\{ f_1,...,f_m \}
_{q(m)} $ of $f_1,...,f_m\in \cal F$ gives only $\bf R$-homogeneous
${\cal A}_p$-additive superdifferential, where $U$ is open in ${\cal
A}_p^n$, a vector $q(m)$ indicates on the order of multiplication.
This superdifferential operator can be extended on the corresponding
family of converging series arising from such final products and
hence on locally analytic functions of $z_1,...,z_n\in {\cal A}_p$.
Since there are many embeddings of $\bf C$ into ${\cal A}_p$, $p\ge
2$, then to encompass the case of complex holomorphic functions in
such theory it was introduced the condition analogous to
holomorphicity: $\partial f/ \partial {\tilde z}=0$, where $z{\tilde
z}=\vert z \vert ^2$, ${\tilde z}$ denotes the adjoint of a
Cayley-Dickson number $z$. \par To make it accurately it was used
the notions of words and phrases and germs of functions and the
Stone-Weierstrass theorem. On the other hand, such formalized
definition of superdifferentiability does not impose from the
beginning the condition of local analyticity. It was proved in
\cite{luoyst,luoystoc} under definite conditions equivalence of
Cayley-Dickson holomorphicity, independence of the line integral
over ${\cal A}_p$ of holomorphic function from the rectifiable paths
with common beginning and final points, local analyticity of
functions. If $f: U\to {\cal A}_p$ is superdifferentiable, where $U$
is open in ${\cal A}_p^n$, then there exists an ${\cal
A}_p$-superdifferentiable $g: U\to {\cal A}_p$ such that $(\partial
g(z)/\partial z).1=f(z)$ and ${\hat f}(z) := (\partial g(z)/\partial
z)$. Then the path integral is defined with the help of the operator
${\hat f}$. The line integral along each rectifiable path $\gamma $
has the continuous extension from $C^1(U,{\cal A}_p)$ onto
$C^0(U,{\cal A}_p)$. Thus there can be given the following.
\par {\bf 1. Definitions.} Define ${\hat f}$ as a generalized function for
each $f\in C^0 (U,{\cal A}_p)$ in the sense of distributions: $(\phi
, {\hat f})_{\gamma } := \int_{\gamma } \phi (z)f(z)dz$, where $\phi
$ is infinite ${\cal A}_p$-superdifferentiable, $\gamma $ is a
rectifiable path in a subset $U$ open in ${\cal A}_p^n$. Denote
${\hat f}^{(0)}(z) := {\hat f}(z)$. Define the $k$-th derivative
$(\phi , {\hat f}^{(k)}.(\mbox{ }^1h,...,\mbox{ }^kh))_{\gamma } :=
(-1)^k (\phi ^{(k)}.(\mbox{ }^1h,...,\mbox{ }^kh),{\hat f})_{\gamma
}$ for each $\gamma $ and each infinite ${\cal
A}_p$-superdifferentiable $\phi $ on $U$ and equal to zero in
neighborhoods of $\gamma (0)$ and $\gamma (1)$, where $\gamma :
[0,1]\to U$ is a rectifiable path, $\mbox{ }^1h,...,\mbox{ }^kh\in
{\cal A}_p$, $0<k\in \bf N$. A net $ \{ {\hat f}^{(k)}_{\alpha }:
{\alpha } \in {\cal Y} \} $ converges to ${\hat f}^{(k)}$ on $U$ in
the sense of distributions, if the net $(\phi ^{(k)}.(\mbox{
}^1h,...,\mbox{ }^kh),{\hat f}_{\alpha })_{\gamma }$ converges for
each infinite ${\cal A}_p$-superdifferentiable function $\phi $ on
$U$ and each rectifiable path $\gamma $ in $U$ and each $\mbox{
}^1h,...,\mbox{ }^kh$, where $0\le k\in \bf Z$, $\phi $ has the
corresponding support, $\cal Y$ is a directed set. The same type of
generalized derivatives of a continuous function $f$ and their
convergence in the sense of distributions we adopt for $f$ instead
of $\hat f$ using the identity $f^{(k)}(z) = {\hat f}^{(k)}(z).1$.
\par This article is devoted
to functions of several Cayley-Dickson variables and investigations
of integral representation formulas for them. Moreover, such
formulas are also obtained for differential forms over ${\cal A}_p$,
$p\ge 2$, where ${\cal A}_2=\bf H$, ${\cal A}_3=\bf O$. There are
well-known integral formulas of the Martinelli-Bochner, Leray,
Koppelman type widely used in complex analysis. Here new generalized
formulas are proved for functions of Cayley-Dickson variables
instead of complex. Moreover, analogs of Stein manifolds over
Cayley-Dickson graded algebras are defined and investigated.
\par The results of this paper it is possible to apply for further
investigations of transformation (super)groups and corresponding to
them (super)algebras of manifolds over Cayley-Dickson algebras as
well as loop spaces, measures and stochastic processes on them,
continuing previous studies of groups of loops and groups of
diffeomorphisms of Riemannian and complex manifolds
\cite{lulgcm,lulsqm,lustptg,lupm}, for they are widely used in
mathematical physics and gauge theories.
\section{Differentiable functions of several Cayley-Dickson variables}
\par {\bf 2.1. Theorem.} {\it Let $U$ be an open subset in ${\cal A}_r$,
$2\le r<\infty $, with a $C^1$-boundary $\partial U$ $U$-homotopic
with a product $\gamma _1\times \gamma _2\times ...\times \gamma
_m$, where $m=2^r-1$, $\gamma _j(\theta )= a_j+\rho _j\exp (2\pi
\theta M_j)$, $M_j\in {\cal I}_r$, $|M_j|=1$, $\theta \in [0,1]$,
$\gamma _j([0,1])\subset U$, $0<\rho _j<\infty $, $j=1,2,...,m$,
where $M_1,...,M_m$ are linearly independent over $\bf R$. Let also
$f: cl (U)\to {\cal A}_r$ be a continuous function on $cl (U)$ such
that $(\partial f(z)/\partial {\tilde z})$ is defined in the sense
of distributions in $U$ is continuous in $U$ and has a continuous
extension on $cl (U)$, where $U$ and $\gamma _j$ for each $j$
satisfy conditions of Theorem 3.9 \cite{luoystoc}, where ${\cal
I}_r:= \{ z\in {\cal A}_r: z + {\tilde z}=0 \} $. Then
$$(1) \quad f(z)=(2\pi )^{-m}\int_{\gamma _m}
(\int_{\gamma _{m-1}}(...(\int_{\gamma _2} (\int_{\gamma _1}f(\zeta
_1).[(\partial _{\zeta _1}Ln(\zeta _1- \zeta _2))M_1^*]).[$$
$$(\partial _{\zeta _2}Ln (\zeta _2-\zeta _3))M_2^*])...
[(\partial _{\zeta _{m-1}}Ln (\zeta _{m-1}-\zeta
_m))M_{m-1}^*])...[\partial _{\zeta _m}Ln (\zeta _m-z))M_m^*$$
$$ - (2\pi )^{-m} \int_U \{ (...((( (\partial {\hat f}(\zeta _1)/
\partial {\tilde \zeta }_1).d{\tilde \zeta }_1)\wedge
 \partial _{\zeta _1}Ln(\zeta _1-\zeta _2))M_1^*)$$
 $$\wedge (\partial _{\zeta _2}Ln (\zeta _2-\zeta _3))M_2^*)\wedge ...)\wedge
(\partial _{\zeta _m} Ln (\zeta _m-z))M_m^* \} , \quad m:=2^p-1. $$}
\par {\bf Proof.} We have the identities
$d_{\zeta } [{\hat f}(\zeta ). \partial _{\zeta } Ln (\zeta -z)]=$
$\{ (\partial {\hat f}(\zeta )/\partial {\tilde \zeta }).d{\tilde \zeta } \}
\wedge \partial _{\zeta }Ln (\zeta -z)$ $+\{
(\partial {\hat f} (\zeta )/\partial {\zeta }).d\zeta \}
\wedge \partial _{\zeta } Ln (\zeta -z) $ and
$d_{\zeta }d_{\zeta }Ln (\zeta -z)|_{\zeta \in \gamma }=0$
for $\zeta $ varying along a path $\gamma $,
where for short $f(z)=f(z,\tilde z)$, since there is the bijection
of $z$ with $\tilde z$ on ${\cal A}_r$. There exists
$\bf R$-homogeneous ${\cal A}_r$-additive operator-valued function
$q(\zeta ,z)$ such that $\partial _{\zeta }Ln (\zeta -z)=
q(\zeta ,z).d\zeta $ (see also \S \S 2.1, 2.2, 2.6 and 2.7
\cite{luoystoc}).
As in \cite{luoyst,luoystoc} ${\hat f}(z,{\tilde z}):=
\partial g(z,{\tilde z})/
\partial z$, where $g(z,{\tilde z})$ is an ${\cal A}_r$-valued
function such that $(\partial g(z,{\tilde z})/\partial z).1=
f(z,{\tilde z})$. Since $\zeta _1$ varies along the path
$\gamma _1$, then $d\zeta _1\wedge d\zeta _1|_{\zeta _1\in \gamma _1}=0$.
Consider $z\in U$ and $\epsilon >0$ such that the torus
${\bf T}(z,\epsilon ,{\cal A}_r)$ is contained in $U$, where
$\partial {\bf T} (z,\epsilon ,{\cal A}_r) =\psi _m\times ... \times \psi _2
\times \psi _1$, $\psi _j$ are of the same form as $\gamma _j$ but with
$z$ instead of $a_j$ and with $\rho _j=\epsilon $.
Applying Stokes formula for regions in ${\bf R}^{2^r}$ and componentwise
to ${\cal A}_r$-valued differential forms we get
$$\int_{\partial U}\omega -\int_{\partial {\bf T}(z,\epsilon ,{\cal A}_r)}
\omega =
\int_{U\setminus {\bf T}(z,\epsilon ,{\cal A}_r)}dw, \mbox{ where}$$
$$w=(...(({\hat f}(\zeta _1).[(\partial _{\zeta _1}
Ln(\zeta _1-\zeta _2))M_1^*]). [(\partial _{\zeta _2}Ln (\zeta
_2-\zeta _3))M_2^*])...). [(\partial _{\zeta _m}Ln (\zeta
_m-z))M_m^*]),$$  $m=2^r-1$. Then $(d Ln \exp (\theta M))M^*=d\theta
$ and ${\hat f}(z).d\theta =f(z)d\theta $, since $M\in {\cal I}_r$,
$M\ne 0$ and $\theta \in \bf R$. In view of Theorems 3.9 and 3.23
\cite{luoystoc} we have that \par $\lim_{\epsilon \to 0, \epsilon
>0} (2\pi )^{-m}\int_{\gamma _m}(...(\int_{\gamma _2} (\int_{\gamma
_1}w))...)= f(z)$ and \par  $\lim_{\epsilon \to 0, \epsilon >0}
\int_{U\setminus {\bf T}(z,\epsilon ,{\cal A}_r)}d\omega
=\int_Ud\omega .$ \\
From this formula $(1)$ follows.
\par {\bf 2.1.1. Remark.} Formula $(2.1)$ is the Cayley-Dickson
algebras' analog of the (complex) Cauchy-Green formula.
Since in the sence of distributions $\partial {\hat f}/\partial {\tilde z}=
\partial (\partial g/\partial {\tilde z})/\partial z$ (see Definition 1),
then from $\partial {\hat f}/\partial {\tilde z}=0$ it follows, that
$\partial {\hat f}.1/\partial {\tilde z}=\partial f/\partial {\tilde z}=0$.
If $\partial f/\partial {\tilde z}=0$, then $g$ can be chosen such that
$\partial g/\partial {\tilde z}=0$ \cite{luoystoc}.
Therefore, from Formula $(2.1)$ it follows, that $f$ is
${\cal A}_r$-holomorphic in $U$ if and only if
$\partial f/\partial {\tilde z}=0$ in $U$.
\par {\bf 2.2.1. Remark.} Instead of curves $\gamma $ of Theorem
$2.1$ above or Theorems $3.9, 3.23$ \cite{luoystoc} it is possible
to consider their natural generalization
$\gamma (\theta )+z_0=z_0+\rho (\theta ) \exp (2\pi S(\theta ))$,
where $\rho (\theta )$ and $S(\theta )$ are continuous functions
of finite total variations, $\theta \in [0,1]\subset \bf R$,
$\rho (\theta )\ge 0$, $S(\theta )\in {\cal I}_r$, $2\le r\le \infty $.
Therefore, $\gamma $ is a rectifiable path.
If $S(0)=S(1)\quad mod({\cal S}_r)$ and $\rho (0)=\rho (1)$, then $\gamma $
is a closed path (loop): $\quad \gamma (0)=\gamma (1)$, where
${\cal S}_r:=\{ z\in {\cal I}_r: |z|=1 \} $,
${\cal I}_r:=\{ z\in {\cal A}_r: z+{\tilde z}=0 \} $.
Consider $S$ absolutely continuous such that there exists
$T\in L^1([0,1],{\cal A}_r)$ for which $S(\theta )=
S(0)+\int_0^{\theta }T(\tau )d\tau $
(see Satz 2 and 3 (Lebesgue) in \S 6.4 \cite{kolmfom})
and let $\rho (\theta )>0$ for each $\theta \in [0,1]$.
Evidently, $Mn:=S(1)-S(0)=\int_0^1T(\tau )d\tau $ is invariant relative
to reparatmetrizations $\phi \in Diff^1_+([0,1])$ of diffeomorphisms
of $[0,1]$ preserving the orientation, $n$ is a real number,
$M\in {\cal S}_r$. Then $\Delta Arg (\gamma )
:=Arg (\gamma )|_0^1= 2\pi \int_0^1T(\tau )d\tau $
(see also Formula $(3.7)$ and \S 3.8.3 \cite{luoystoc}).
In view of Theorem $3.8.3$ \cite{luoystoc} for each loop
$\gamma :$ $\quad \Delta Arg (\gamma )\in {\bf Z} {\cal S}_r$.
For each $\epsilon >0$ for the total variation there is
the equality $V(\gamma \epsilon )=V(\gamma )\epsilon $.
Since $\gamma ([0,1])$ is a compact subset in ${\cal A}_r$,
then there exists $\rho _m:=\sup _{\theta \in [0,1]}\rho (\theta )<\infty $.
Hence $z_0+(\gamma \epsilon )([0,1])\subset
B({\cal A}_r,z_0,\rho _m\epsilon )$.
\par Therefore, Theorems $3.9, 3.23, 3.28$ and Formulas
$(3.9, 3.34.i)$ \cite{luoystoc} and Theorem $2.1$ above are true for such
paths $\gamma $ also and Formula $(3.9)$ \cite{luoystoc} takes the form
$$(1)\quad f(z)M=(2\pi n)^{-1}(\int_{\psi }f(\zeta )
(\zeta -z)^{-1}d\zeta ),$$
where $0\ne n\in \bf Z$ for a closed path $\gamma $, $M\in {\cal S}_r$,
Formula $(1)$ generalizes Formula $(3.9)$, when $|n|>1$.
When ${\hat I}n(0,\gamma )=0$, then
$(\int_{\psi }f(\zeta )(\zeta -z)^{-1}d\zeta )=0$
(see also \S 3.23 \cite{luoystoc}).
\par {\bf 2.2.2. Note and Definition.} Let $\Lambda $
denotes a Hausdorff topological space with nonnegative measure $\mu
$ on a $\sigma $-algebra of all Borel subsets such that for each
point $x\in \Lambda $ there exists an open neighborhood $U\ni x$
with $0<\mu (U)<\infty $. Consider a set of generators with real
algebra $\{ i_x: x\in \Lambda \} $ such that $i_xi_y=-i_yi_x$ for
each $x\ne y\in \Lambda \setminus \{ 0 \} $ and $i_x^2=-1$ for each
$x\in \Lambda \setminus \{ 0 \} $, where $0$ is a marked point in
$\Lambda $. Add to this set the unit $1=:i_0$ such that $ai_x=i_xa$
for each $a\in \bf R$ and $x\in \Lambda $. In the case of a finite
set $\Lambda $ the Cayley-Dickson algebra generated by such
generators is isomorphic with ${\cal A}_{N-1}$, where $N=card
(\Lambda )$ is the cardinality of the set $\Lambda $. \par For the
infinite subset of generators  $ \{ i_0, i_{x_j}: j\in {\bf N},
x_j\in \Lambda \} $ the construction from \S 3.6.1 \cite{luoystoc}
produces the algebra isomorphic with ${\cal A}_{\infty }$.
Therefore, consider the case $card (\Lambda )>\aleph _0$. Due to the
Kuratowski-Zorn lemma we can suppose, that $\Lambda $ is linearly
ordered and this linear ordering gives intervals $(a,b):=\{ x\in
\Lambda : a<x<b \} $ being $\mu $-measurable, for example, $\Lambda
= {\bf R}^n\times ({\bf R}/{\bf Z})^m$ has the natural linear
ordering induced by the linear ordering from $\bf R$ and by the
lexicographic ordering in the product, where $n, m \in \bf N$.
\par Then consider a finite partition $\Lambda $ into a disjoint union
$\Lambda =\bigcup_{j=0}^pA_j$, where $x<y$ for each $x\in A_j$ and
$y\in A_l$ for $j<l\le p$, $p\in \bf N$, $0\in A_0$. The family of
such partitions we denote $\cal Z$. Let $T\in \cal Z$, $x_j\in A_j$
be marked points. Then there exists a step function $f_T$ such that
$f_T(x)=C_ji_{x_j}$ for each $x\in A_j$, where $C_j\in \bf R$.
Consider the norm $\| f_T \|_{\Lambda } ^2:=\int_{\Lambda }
f_T(x){\tilde f}_T(x)\mu (dx)$, where ${\tilde f}_T(x) := C_0\chi
_{A_0}(x)\delta _{0,x_0} - \sum_{x_j\ne 0}C_j\chi _{A_j}(x)i_{x_j}$,
$\chi _A(x)=1$ for $x\in A$, $\chi _A(x)=0$ for $x\notin A$, $\delta
_{x,y}=1$ for $x=y$, $\delta _{x,y}=0$ while $x\ne y$. To each $f_T$
put the element $z_{f_T}:=\sum_j C_ji_{x_j}\mu (A_j)$. \par The
algebra which is the completion by the norm $\| * \|_{\Lambda } $ of
the minimal algebra generated by the family of elements $z_{f_T}$
for $f_T$ from the family $\cal F$ of all step functions and all
their ordered final products we denote by ${\cal A}_{\Lambda }$.
\par {\bf 2.2.3. Theorem.} {\it The set ${\cal A}_{\Lambda }$
is the power-associative noncommutative nonassociative algebra over
$\bf R$ complete relative to the norm $\| * \|_{\Lambda }$ with the
centre $Z({\cal A}_{\Lambda })=\bf R$, moreover, there are
embeddings ${\cal A}_{\infty }\hookrightarrow {\cal A}_{\Lambda }$
for $card (\Lambda )\ge \aleph _0$. The set of generators of the
algebra ${\cal A}_{\Lambda }$ has the cardinality $card (\Lambda )$
for $card (\Lambda )\ge card ({\bf N})$. There exists the function
${\overrightarrow { \exp }}(\int_{\Lambda }f(x)\mu (dx))$ of the
ordered integral product from ${\cal A}_{\Lambda }$ onto ${\cal
A}_{\Lambda }$.}
\par {\bf Proof.} For $card (\Lambda )\le \aleph _0$ the algebra
${\cal A}_{\Lambda }$ is  isomorphic with ${\cal A}_{N-1}$ or ${\cal
A}_{\infty }$. Thus it remains to consider the case $card (\Lambda
)>\aleph _0$. For each $f_T\in \cal F$ it can be defined the ordered
integral exponential product \\
${\overrightarrow {\exp }} (\int_{\Lambda }f_T(x)\mu (dx)) := \{
\exp (C_0\mu (A_0)\pi i_{x_0}/2)...\exp (C_p\mu (A_p) \pi i_{x_p}/2)
\}_{q(p+1)} $ \\
with $q(p+1)$ corresponding to the left order of brackets. Thus
there exist the embeddings of ${\cal A}_{\infty }$ into ${\cal
A}_{\Lambda }$. Then $Z({\cal A}_{\Lambda })=\bf R$. The completion
of the family ${\cal F}$ contains all functions of the type $f(x) =
\sum_j f_j(x)\chi _{A_j}(x) i_{x_j}$, where $\{ A_j: j\in {\bf N} \}
$ is the disjoint union of $\Lambda $, each $A_j$ is $\mu
$-measurable, $f_j\in L^2(\Lambda ,\mu ,{\bf R})$ and $\lim_{n\to
\infty } \sum_{j>n} \| f_j(x) \chi _{A_j}(x)\|^2_{L^2(\Lambda , \mu
, {\bf R})}=0$. \par  Since $\exp (M)= \cos (|M|) + M \sin (|M|)
/|M|$ for each $M\in {\cal A}_{\infty }$ and $|\exp (M)-1|\le \exp
(|M|)-1$, then for each $f\in {\cal A}_{\Lambda }$ there exists
\par $\lim_{{\cal F}\ni f_T\to f} {\overrightarrow {\exp }}($
$\int_{\Lambda }f_T(x)\mu (dx))=: {\overrightarrow {\exp
}}(\int_{\Lambda }f(x)\mu (dx))$ \\
relative to $\| * \|_{\Lambda }$.
From $\exp (\pi i_x/2)=i_x$ for each $x\in \Lambda \setminus \{ 0 \}
$ it follows, that the family of all elements of the type
${\overrightarrow { \exp }}(\int_{\Lambda }f(x)\mu (dx))$, $f\in
\cal F$ contains all generators of the embedded subalgebra ${\cal
A}_{\infty }$ generated by the countable subfamily $\{ i_{x_j}: j\in
{\bf N} \} $. \par  The completion ${\tilde {\cal F}}$ of the family
$\cal F$ by the norm $\| * \|_{\Lambda }$ is the infinite
dimensional linear subspace over $\bf R$ in ${\cal A}_{\Lambda }$.
All possible final ordered products from ${\tilde {\cal F}}$ and the
completion of their $\bf R$ -linear span by the norm $\|
* \|_{\Lambda }$ produces ${\cal A}_{\Lambda }$. Then for each
element from ${\cal A}_{\Lambda }$ there exists the representation
in the form of the ordered integral exponential product. Since
${\cal A}_{\Lambda }$ is the algebra over $\bf R$ and $card (\Lambda
)^{\aleph _0}=card (\Lambda )$, then the family of generators of the
algebra ${\cal A}_{\Lambda }$ has the cardinality $card (\Lambda )$.
\par {\bf 2.2.4. Note.} Evidently Propositions $2.2.1, 2.3, 2.6$ and Corollary
$2.4$, Lemma $2.5.1$ \cite{luoystoc} are accomplished in the case of
${\cal A}_{\Lambda }$ with ${\bf b}={\bf b}_{\Lambda }$ instead of
${\bf b}={\bf b}_r$. Definition $2.5$ has the meaning also for
${\cal A}_{\Lambda }$. Theorem $2.7$ is also accomplished for ${\cal
A}_{\Lambda }$, since for each $z\in {\cal A}_{\Lambda }$ there
exists the embedded subalgebra isomorphic with ${\cal A}_{\infty }$
and containing $z$. A path $\gamma $ is rectifiable, hence it has a
countable dense subset. For each $\epsilon
>0$ there exists a subalgebra isomorphic with ${\cal A}_{\infty }$ the projection
$\psi (t)$ on which of the path $\gamma $ differs from $\gamma (t)$
no more than on $\epsilon $ for each $t\in [a,b]$, where $\gamma :
[a,b]\to {\cal A}_{\Lambda }$. For ${\cal A}_{\infty }$ with the
help of the projections $P_r$ we have $\psi =\lim_{r\to \infty
}P_r(\psi )$, $P_r(\psi )\subset U_r$, $ \{ P_r(\gamma ): r \in {\bf
N} \} $, converges to $\psi $ uniformly on the compact segment
$[a,b]\subset \bf R$, where $U_r=P_r(U)$. Take a sequence of such
path $\psi _n$ with $\sup_{t\in [a,b]} |\psi _n(t)-\gamma (t)|<1/n$.
Then $\int_{\psi _n}f(z)dz$ forms the Cauchy sequence in ${\cal
A}_{\Lambda }$, which is complete. Thus there exists $\lim_{n\to
\infty }\int_{\psi _n}f(z)dz=\int_{\gamma }f(z)dz$. Consequently,
the integral along the path has the unique continuous extension on
$C^0_b(U,{\cal A}_{\Lambda })$. For a continuous function $f$ on an
open domain $U$ in ${\cal A}_{\Lambda }$ there exists a generalized
operator $\hat f$ in the sense of distributions on rectifiable paths
in $U$. Then Definition 1 has the natural extension on ${\cal
A}_{\Lambda }$.
\par In Note $2.8$ it can be used $l_2({\bf R})^m$ instead of
$\bf R^{2^rm}$ and represent the differential forms $\eta $ over
${\cal A}_{\infty }$ as the pointwise limits (or to use the uniform
converegence on compact subsets) of differential forms over ${\cal
A}_r$ for $r$ tending to the infinity, since $z_r\to z$ while $r$
tends to the infinity, where $z\in {\cal A}_{\infty }$,
$z_r:=P_r(z)$. In the case of ${\cal A}_{\Lambda }$ when $card
(\Lambda )> \aleph _0$ this can be used pointwise, since for each
$z\in {\cal A}_{\Lambda }$ there exists a subalgebra isomorphic with
${\cal A}_{\infty }$ and containing $z$. In the general case: \\
$(i)\quad \eta (z,{\tilde z})=\sum_{I,J} \eta _{I,J} \{ (d\mbox{
}^{p_1}z^{\wedge I_1}\alpha _1 \wedge ... \wedge d\mbox{
}^{p_n}z^{\wedge I_n}\alpha _n \wedge d\mbox{ }^{t_1}{\tilde
z}^{\wedge J_1}\beta _1\wedge ... \wedge
d\mbox{ }^{t_n}{\tilde z}^{\wedge J_n}\beta _n \} _{q(|I|+|J|+2n)}$ \\
this is the differential form over ${\cal A}_{\Lambda }$, where each
$\eta _{I,J}(z,{\tilde z})$ is a continuous function on an open
subset $U_n$ in ${\cal A}_{\Lambda }^n$ with values in ${\cal
A}_{\Lambda }$, $I=(I_1,...,I_n)$, $J=(J_1,...,J_n)$,
$|I|:=I_1+...+I_n$, $1\le p_1\le p_2\le ... \le p_n\in \bf N$, $1\le
t_1\le t_2\le ... \le t_n\in \bf N$, $0\le I_k \in \bf Z$, $0\le
J_k\in \bf Z$, $\alpha _k, \beta _k\in {\cal A}_{\Lambda }$ are
constants for each $k=1,...,n$, $d\mbox{ }^pz^0:=1$, $d\mbox{
}^p{\tilde z}^0:=1$, $n\in \bf N$, $\pi ^l_n(U_l)\subset U_n$ for
all $l\ge n$, where $\pi ^l_n: {\cal A}_{\Lambda }^l\to {\cal
A}_{\Lambda }^n$ is the natural projection for each $l\ge n$. The
converegence on the right side of Formula $(i)$ in the case of an
infinite series by $I$ or $J$ is supposed relative to $C^0_b(W,{\cal
A}_{\Lambda }^{\wedge *})$-topology of the uniform convergence on
$W$, where $W=pr-\lim \{ U_n, \pi ^l_n, {\bf N} \} $, ${\cal
A}_{\Lambda }^{\wedge *}$ is supplied with the norm topology
inherited from the topologically adjoint space of all poly $\bf
R$-homogeneous ${\cal A}_{\Lambda }$-additive functionals.
\par In Note $2.10$ \cite{luoystoc} define ${\cal A}_{\Lambda ,s,p}$ an use
projections $\pi _{s,p,t}$ for each $s\ne p\in \bf b$. Theorems
$2.11, 2.15$ and Corollaries $2.13, 2.15.1$ are transferrable on
${\cal A}_{\Lambda }$ with $card (\Lambda )\ge \aleph _0$ by
imposing the condition of $(2^r-1)$-connectedness $P_r(U)=:U_r$ for
each $r\ge 3$ and every embedding of ${\cal A}_{\infty }$ into
${\cal A}_{\Lambda }$ for $card (\Lambda )>\aleph _0$ while
corresponding ${\cal A}_r\subset {\cal A}_{\infty }$ considering
$\pi _{s,p,t}(U)$ for each $s=i_{2k}$, $p=i_{2k+1}$, $0\le k\in \bf
Z$. Then Definitions $2.12, 2.14$ and Theorem $2.16$, Notes $2.17,
3.1$ are also accomplished for ${\cal A}_{\Lambda }$. Corollary
$3.3$ in this case follows from Theorem $3.6.2$ \cite{luoystoc}.
\par As for ${\cal A}_{\infty }$ the algebra ${\cal A}_{\Lambda }$
with $card (\Lambda )> \aleph _0$ has not any finite and even
countable set of constants $\{ a_s, b_s \} $ in ${\cal A}_{\Lambda
}$ such that $z^*=\sum_s a_szb_s$ for each $z\in {\cal A}_{\Lambda
}$ could be written in as such series or sum. That is the algebraic
antiautomorphism of order two $\theta (z):=z^*$ with $\theta \circ
\theta =id$ is not internal in ${\cal A}_{\Lambda }$ and indeed $z$
and $z^*$ are algebraically independent variables in such infinite
dimensional Cayley-Dickson algebra.
\par {\bf 2.2.5. Proposition.} {\it Let $U$ be an open subset in
${\cal A}_p$ and $f: U\to {\cal A}_p$ be a function on $U$, where
$2\le p\le \infty $, let also ${\cal A}_{\Lambda }$ be the
Cayley-Dickson algebra as in \S 2.2.2 with $card (\Lambda )\ge
\aleph _0$, then $f$ is $z$-superdifferentiable if and only if there
exists an open subset $W$ in ${\cal A}_{\Lambda }$ and a
$z$-superdifferentiable function $g: W\to {\cal A}_{\Lambda }$ such
that its restriction on $U$ coincides with $f$, $g|_U=f$. }
\par {\bf Proof.} In view of Theorem 2.2.3 there exists the embedding of
${\cal A}_p$ into ${\cal A}_{\Lambda }$. If $g$ is
$z$-superdifferentiable on $W$, then from the definition it follows,
that $g|_V$ is $z$-superdifferentiable, where $V=W\cap {\cal A}_p$
and $V$ is open in ${\cal A}_p$. Vice versa if $f$ is
$z$-superdifferentiable on $U$, then it is locally $z$-analytic on
$U$ (see Theorems 2.15 and 3.10 \cite{luoystoc}). For each $z_0\in
U$ there exists a power series in $(z-z_0)$ converging in a ball
$B({\cal A}_p,z_0,r)$ with the centre $z_0$ and positive radius
$r>0$ the expansion coefficients of which belong to ${\cal A}_p$.
Therefore, in the variable $(z-z_0)$ this series uniformly converges
in $B({\cal A}_{\Lambda },z_0,r')$ for each $0<r'<r$. The union of
such balls $B({\cal A}_{\Lambda },z_0,r')$ is the open subset in
${\cal A}_{\Lambda }$ which we denote by $W$. On each open
intersection of each corresponding pair of such balls the functions
given by such series coincide, that gives the
$z$-superdifferentiable function $g$ on $W$ with $g|_U=f$. Certainly
there can be found others $z$-superdifferentiable extensions $g$ of
$f$.
\par {\bf 2.3. Theorem.} {\it Let $U$ be a bounded open subset
in ${\cal A}_{\Lambda }$ and let $f: U\to {\cal A}_{\Lambda }$ be a
bounded continuous function. Then there exists a continuous function
$u(z)$ which is a solution of the equation
$$(1)\quad (\partial u(z)/\partial {\tilde z})={\hat f}$$
in $U$, in particular, $(\partial u(z)/\partial {\tilde
z}).1=f(z)$.}
\par {\bf Proof.} Using embeddings of ${\cal A}_p$, $2\le p\in \bf N$, into
${\cal A}_{\Lambda }$ it is sufficient to prove this statement for
arbitrary $2\le p\in \bf N$. Take $2\le p\in \bf N$ and ${\cal A}_p$
one-forms $d\zeta _l$ expressible through $d\zeta $ as
$\sum_{j=1}^{k(l)}P_{j,1,l}d\zeta P_{j,2,l}$ with fixed nonzero
$P_{j,q,l}\in {\cal A}_p$, where $l=1,2,3,...,2^p$, $k(l)\in \bf N$.
Choose $d\zeta _l$ to be satisfying conditions $d\zeta _{2^p}\wedge
\nu =\xi (z)(...((d\mbox{ }^1z\wedge d\mbox{ }^2z )\wedge d\mbox{
}^3z)...) \wedge d\mbox{ }^{2^p}z)$, $d{\tilde {\zeta }}_{2^p}\wedge
\nu =0$, where
$$\nu =(\partial _{\zeta _1} Ln(\zeta _1-\zeta
_2))M_1^*]). [(\partial _{\zeta _2}Ln (\zeta _2-\zeta
_3))M_2^*])...). [(\partial _{\zeta _m}Ln (\zeta _m-z))M_m^*]),$$
$m=2^p-1$ as in \S 2.1, $z=\sum_{l=1}^{2^p}\mbox{ }^lzi_{l-1}\in U$,
$\mbox{ }^lz\in \bf R$, $i_0,...,i_{2^{p-1}}$ are generators of
${\cal A}_p$, $\xi : U\to {\hat {\cal A}}_p$ is a function nonzero
and finite almost everywhere on $U$ relative to the Lebesgue
measure. Then there exist $d\zeta _{2^p}$ and $\nu $ such that the
continuous function
$$(2)\quad u(z):=-(2\pi )^{1-2^p}\int_U ({\hat f}
(\zeta _1).d{\tilde \zeta }_{2^p}) \wedge \nu $$ is a solution of
equation $(1)$. To demonstrate this take closed curves (paths)
$\gamma _j$ in $U$ as in \S 2.1 and \S 2.2.4, for example, such that
$\zeta _j\in \gamma _j$ satisfy conditions: $(\zeta _{2s-1}-\zeta
_{2s})=\eta _{2s-1}$ for each $s=1,...,2^{p-1}-1$, $\zeta _{2^p-1}=-
r\eta _{2^p-1}^*$ with $0<r<1$, where $\eta _{2s-1}$ and $\eta
_{2s}$ belong to the plane $i_{2s-1}{\bf R}\oplus i_{2s}{\bf R}$,
$\eta _{2s}={\tilde \eta }_{2^s-1}$ for each $s=1,...,2^{p-1}$,
${\tilde {\eta _1}}={\eta _{2^p}}$. Hence $d\eta _1\wedge d\eta
_1=0$, $d\eta _2\wedge d\eta _3^*=0$, ...,$d\eta _{2^p-1}\wedge
d\eta _{2^p}^*=0$, $d\eta _{2^p}\wedge d\eta _{2^p}=0$,
\par $(i)$ $\eta _v^kd\eta _v=(d\eta _v)\eta _v^k$ for each $v$ \\
for $k=1$ and $k=-1$. \par These variables are expressible as $\zeta
_l=\sum_{j=1}^{k(l)}P_{j,1,l}\zeta P_{j,2,l}$ (see \S \S 3.7 and
3.28 \cite{luoyst}). Therefore, there exists a subgroup of the group
of all ${\cal A}_p$-holomorphic diffeomorphisms of $U$ preserving
Conditions $(ii)$ and the construction given above has natural
generalizations with the help of such diffeomorphisms.
\par Supoose at first, that $f$ is continuously differentiable in $U$. Each
$\zeta _j$ is expressible in the form $\zeta _j= \sum_l\mbox{
}^lb_jS_l$, where $\mbox{ }^lb_j\in \bf R$ are real variables,
$S_l\in \{ i_0,i_1,...,i_{2^p-1} \} $, hence differentials
$(\partial f/\partial \zeta _j).d\zeta _j=\sum_l \{ (\partial
f/\partial z).S_ld\mbox{ }^lb_j$ $+(\partial f/\partial {\tilde z}).
{\tilde S_l}d\mbox{ }^lb_j \} $ are defined. Consider a fixed
$z_0\in U$. We take a $C^{\infty }$-function $\chi $ on ${\cal A}_p$
such that $\chi =1$ in a neighbourhood $V$ of $z_0$, $V\subset U$,
$\chi =0$ in a neighbourhood of ${\cal A}_p\setminus U$. Then
$u=u_1+u_2,$ where
$$u_1(z):=-(2\pi )^{1-2^p}\int_U [\chi (\zeta _1)
({\hat f}(\zeta _1).d{\tilde \zeta }_{2^p})] \wedge \nu ,$$
$$u_2(z):=-(2\pi )^{1-2^p}\int_U [(1-\chi (\zeta _1))
({\hat f}(\zeta _1).d{\tilde \zeta }_{2^p})]\wedge \nu .$$ Then
$$u(z):=-(2\pi )^{1-2^p}\int_{{\cal A}_p}[\chi (\zeta _1+z)
({\hat f}(\zeta _1+z).d{\tilde \zeta }_{2^p})]\wedge \psi ,\mbox{
where}$$
$$\psi :=(\partial _{\zeta _1} Ln(\zeta _1-\zeta _2))M_1^*]). [(\partial
_{\zeta _2}Ln (\zeta _2-\zeta _3))M_2^*])...). [(\partial _{\zeta
_m}Ln (\zeta _m))M_m^*]).$$ Since $\partial _{\zeta _{2^p}} \{ [\chi
(\zeta _1+z){\hat f}(\zeta _1+z).]\wedge \psi \} =0$ and $(\partial
/\partial {\tilde \zeta _{2^p}}) \{ [\chi (\zeta _1+z){\hat f}(\zeta
_1+z).]\wedge \psi \}.d{\tilde \zeta }_{2^p}$ $=\partial _{{\tilde
\zeta }_{2^p}} \{ [\chi (\zeta _1+z){\hat f}(\zeta _1+z).]\wedge
\psi \},$ then due to Equations $(i,ii)$ \\
$(\partial u(z)/\partial {\tilde z})= -(2\pi )^{1-2^p}\int_{{\cal
A}_p}\partial _{{\tilde \zeta }_{2^p}}
\{ [\chi (\zeta _1+z){\hat f}(\zeta _1+z)]\wedge \psi \} $.  \\
In view of Theorem 2.1 applied to ${\hat f}.S$ for each $S\in \{
i_0,i_1,...,i_{2^p-1} \} $ we have $(\partial u_1/\partial {\tilde
z})=\hat f$ in $V$, consequently, $(\partial u/\partial {\tilde
z})=\hat f$ in a neighbourhood of $z_0$.
\par Taking a sequence
$f^n$ of continuously differentiable functions uniformly converging
to $f$ on $U$ we get the corresponding $u^n$ such that in the sence
of distributions $(\partial u/\partial {\tilde z})= \lim_{n\to
\infty } (\partial u^n/\partial {\tilde z}) =\lim_n{\hat f}^n=\hat
f$. Consider the family of all embeddings $\theta $ of ${\cal A}_p$
to ${\cal A}_{\Lambda }$, $2\le p<\infty $. There exists the
generalized function (operator) $\hat f$ on $U$, hence there exists
the restriction ${\hat f}_{p,\theta }$ on $U\cap \theta ({\cal
A}_p)$ for each $(p,\theta )$. As the function this restriction
evidently exists. For distributions it is possible to take them on a
base space of cylindrical functions on the algebra of cylindrical
subsets with bases in the projection $\theta ({\cal A}_p)$. Each
rectifiable path $\gamma $ is the limit of the uniformly converging
net of paths $\gamma _{p,\theta }$, since $\gamma ([0,1])$ is
compact. Therefore, such restriction exists in the sence of
distributions.
\par Thus there exists the solution $u_{p,\theta }$
of $(1)$ given by $(2)$ on $U\cap \theta ({\cal A}_p)$. The family
of all $(p, \theta )$ is directed: $(p_1, \theta _1)\le (p_2, \theta
_2)$ if and only if $p_1\le p_2$ and $\theta _1({\cal
A}_{p_1})\subset \theta _2({\cal A}_{p_2})$. Since ${\hat
f}_{p,\theta }$ converges to $\hat f$ in the sense of distributions
by the ultrafilter of the set $\{ (p,\theta ) \} $, then
$u_{p,\theta }$ converges to the solution $u$ on $U$, since there
exists $\partial u(z)/\partial {\tilde z}={\hat f}(z)$.
\par {\bf 2.4. Theorem.} {\it Let $U$ be an open subset in ${\cal A}_p^n$,
$2\le p\in \bf N$, $n\in \bf N$. Then for every compact subset $K$
in $U$ and every multi-order $k=(k_1,...,k_n)$, there exists a
constant $C>0$ such that
$$\max_{z\in K} |\partial ^kf(z)|\le C\int_U|f(z)|d\sigma _{2^pn}$$
for each ${\cal A}_p$-holomorphic function $f$, where $d\sigma
_{2^pn}$ is the Lebesgue measure in ${\cal A}_p^n$.}
\par {\bf 2.5. Corollary.} {\it Let $U$ be an open subset in ${\cal
A}_p^n$, $2\le p\in \bf N$, $n\in \bf N$, and let $f_l$ be a
sequence of ${\cal A}_p$-holomorphic functions in $U$ which is
uniformly bounded on every compact subset of $U$. Then there is a
subsequence $f_{k_j}$ converging uniformly on every compact subset
of $U$ to a limit in $C^{\omega }_z(U,{\cal A}_p)$.}
\par Proofs of Theorem 2.4 and Corollary 2.5 follow from
Theorem 2.1 above and Theorem 3.9 \cite{luoyst} (see also
\cite{luoystoc}) analogously to Theorem 1.1.13 and Corollary 1.1.14
\cite{henlei}.
\par {\bf 2.6. Definitions.} Let $U$ be an open subset in ${\cal A}_p^n$
and $f: U\to {\cal A}_p^m$ be an ${\cal A}_p$-holomorphic function,
then the matrix: $J_f(z):=(\partial f_j(z)/\partial z_k)$ is called
the ${\cal A}_p$-Jacobi matrix, where $j=1,...,m$, $k=1,...,n$. To
this operator matrix there corresponds a real $(2^pm)\times
(2^pn)$-matrix while $2\le p\in \bf N$ or operator from $X^m$ into
$X^n$ of the underlying real Hilbert space $X$ of ${\cal A}_p$ for
infinite $p=\Lambda $. Denote by $rank_{\bf R}(J_f(z))$ a rank of a
real matrix or operator corresponding to $J_f(z)$. This rank may be
infinite. Then $f$ is called regular at $z\in U$, if $rank_{\bf
R}(J_f(z))=2^p\min (n,m)$ for finite $p$ or $ker (f'(z))=0$ and
$Range (f'(z))$ is algebraically isomorphic with ${\cal A}_{\Lambda
}^m$ such that $Range (f'(z))\oplus {\cal A}_{\Lambda }^{n-m}= {\cal
A}_{\Lambda }^n$ when $m\le n$ or $Range (f'(z))=X^n$ and $ker
(f'(z))$ is algebraically isomorphic with ${\cal A}_{\Lambda
}^{m-n}$ while $m>n$. If $U$ and $V$ are two open subsets in ${\cal
A}_p^n$, then a bijective surjective mapping $f: U\to V$ is called
${\cal A}_p$-biholomorphic if $f$ and $f^{-1}: V\to U$ are ${\cal
A}_p$-holomorphic.
\par {\bf 2.7. Proposition.} {\it Let $U$
and $V$ be open subsets in ${\cal A}_p^n$ and ${\cal A}_p^m$
respectively. If $f: U\to {\cal A}_p^m$ and $g: V\to {\cal A}_p^k$
are ${\cal A}_p$-holomorphic functions such that $f(U)\subset V$,
then $g\circ f: U\to {\cal A}_p^k$ is ${\cal A}_p$-holomorphic and
$J_{g\circ f}(z)=J_g(f(z)).(J_f(z).h)$ for each $h\in {\cal
A}_p^n$.}
\par {\bf Proof.} In view of Definition 2.2 and Theorems 2.15 and 3.10
\cite{luoystoc} $(\partial g_j(f(z))/\partial z_l).\zeta =
\sum_{s=1}^m\sum_{l=1}^k (\partial g_j(\xi )/ \partial \xi _s)|_{\xi
=f(z)}. (\partial f_s(z)/\partial z_l).h_l$, where
$h=(h_1,...,h_n)$, $h_l\in {\cal A}_p$ for each $l=1,...,n$, since
$f(U)\subset V$ and this is evident for ${\cal A}_p$-polynomial
functions and hence for locally converging series of ${\cal
A}_p$-holomorphic functions.
\par {\bf 2.8. Proposition.} {\it Let $U$ be a neighbourhood of
$z\in {\cal A}_p^n$ and let $f: U\to {\cal A}_p^n$ be an ${\cal
A}_p$-holomorphic function. Then $f$ is ${\cal A}_p$-biholomorphic
in some neighbourhood $W$ of $z$ if and only if $f$ is regular at a
point $z\in U$.}
\par {\bf Proof.} From Proposition 2.7 it follows, that
the condition of regularity of $f$ on $U$ is necessary. Prove the
sufficiency. In view of Definition 2.2, Theorems 2.15 and 3.10 and
Note 3.11 \cite{luoystoc} an incerement of $f$ can be written in the
form $f(z+\zeta )=f(z)+f'(z).\zeta + O(|\zeta |^2)$ for each $\zeta
\in {\cal A}_p^n$ such that $z+\zeta \in U$. Then there exists a
neighborhood $W\supset B(z,2\epsilon ,{\cal A}_p^n)$ in which
$|g(z+\zeta )|\le C|\zeta |^2$, where $0<\epsilon <(2C)^{-1}$, $C$
is a positive constant, $g:=id-f$. Thus there exists an ${\cal
A}_p$-holomorphic function $w$ on an open neighbourhood $W$ of $z$
in $U$ such that $w$ is given by the series $w=\sum_{k=1}^{\infty
}g_k$, where $g_{k+1}=g\circ g_k$ for each $k\in \bf N$ and
$g_1:=g$, $g:=id-f$, since for each $\eta \in W$ there exists $r>0$
such that $B(\eta ,r,{\cal A}_p)\subset U$ and the series for $w$ is
convergent on $B(\eta ,r,{\cal A}_p)$ with $w (B(z,\epsilon ,{\cal
A}_p^n))\subset B(z,2\epsilon ,{\cal A}_p^n)$.
\par Since $f'(z)$
is the continuous epimorphism from ${\cal A}_p^n$ onto ${\cal
A}_p^n$, then its graph is closed. On the other hand, $f'(z)$ is
bijective and there exists the $\bf R$-linear operator
$(f'(z))^{-1}$. The graph of it $Gr (f'(z))^{-1} = \{ (x,y):
x=f'(z).y; x, y \in {\cal A}_p^n \} $ is closed in ${\cal
A}_p^n\otimes {\cal A}_p^n$, since the graph of $f'(z)$ is closed.
In view of the closed mapping theorem (see 14.3.4 \cite{nari})
$(f'(z))^{-1}$ is continuous. Thus the operator $f'(z)$ is
invertible. In view of the inverse mapping theorem (see \S X.7
\cite{zorich}) there exists $f^{-1}$ continuusly (Frech\'et)
differentiable on a neighborhood $W$ of $f(z)$. Since $\partial
f(z)/\partial {\tilde z}=0$, then $\partial f^{-1}(\zeta )/\partial
{\tilde {\zeta }}=0$ on $W$.
\par For $\eta $ in a sufficiently
small neighborhood $W$ of $z$ there is satisfied the inequality $\|
1 - f'(z)^{-1} f'(\eta ) \| < 1$, consequently, $f'(\eta )$ is
invertible for each $\eta \in W$. The operator $f'(\eta )$ is
continuous by $\eta $ on $U$, hence there exists a neighbourhood $V$
of $z$ such that $f$ is regular on $V$, since $f'(\eta )$ is $\bf
R$-homogeneous and ${\cal A}_p$-additive and $f'(z)({\cal
A}_p^n)={\cal A}_p^n$. Hence $f(V)$ is open in ${\cal A}_p^n$. Since
$w$ is the limit of the uniformly convergent series of ${\cal
A}_p$-holomorphic functions, then $w$ is ${\cal A}_p$-holomorphic on
$W$. From $(id+h)\circ f=f\circ (id+h)=id$ on $B(z,\epsilon ,{\cal
A}_p^n)$ it follows, that $f$ is ${\cal A}_p$-biholomorphic on a
neighbourhood of $z$.
\par {\bf 2.9. Corollary.} {\it Let $X$ be a subset in ${\cal
A}_p^n$, $2\le p<\infty $ or $p=\Lambda $, and $k\in \{ 1,2,...,n-1
\} $, then the following conditions are equivalent:
\par $(i)$ for each $\zeta \in X$ there exists an ${\cal A}_p$-biholomorphic
map $f=(f_1,...,f_n)$ in some neighbourhood $U$ of $\zeta $ such
that $f$ is regular on $U$ and $X\cap U= \{ z\in U:
f_{k+1}(z)=0,...,f_n(z)=0 \} $;
\par $(ii)$ for each $\zeta \in X$ there exists a neighbourhood
$V$ of $\zeta $ and a regular ${\cal A}_p$-holomorphic map $g: V\to
{\cal A}_p^{n-k}$ such that $X\cap V= \{ z\in V: g(z)=0 \} $.}
\par {\bf Proof.} The implication $(i)\Rightarrow (ii)$ follows by
taking $g=(f_{k+1},...,f_n)$ on $V=U$. To prove implication
$(ii)\Rightarrow (i)$ take $\zeta \in X$, $g$ and $V$ as in $(ii)$.
There exists the $\bf R$-linear operator $G'$ corresponding to
$g'(\zeta )$ from ${\cal A}_p^n$ onto ${\cal A}_p^{n-k}$. Thus there
exists a right ${\cal A}_p$-superlinear operator $P$ from ${\cal
A}_p^n$ onto ${\cal A}_p^k$ such that $P\oplus g'(\zeta )$ from
${\cal A}_p^n$ onto ${\cal A}_p^n$ is invertible. The graph of it
$Gr ((P\oplus g'(\zeta ))^{-1} = \{ (x,y): x=(P\oplus g'(\zeta
))^{-1}.y; x, y \in {\cal A}_p^n \} $ is closed in ${\cal
A}_p^n\otimes {\cal A}_p^n$, since the graph of $(P\oplus g'(\zeta
))$ is closed. In view of the closed mapping theorem (see 14.3.4
\cite{nari}) $(P\oplus g'(\zeta ))^{-1}$ is continuous. Thus the
operator $(P\oplus g'(\zeta ))$ is invertible. In view of the
implicit mapping theorem and addition 3 to it (see \S X.7
\cite{zorich}) there exists $(P\oplus g)^{-1}$ continuusly
(Frech\'et) differentiable on a neighborhood $W$ of $(P\oplus
g)(\zeta )$.
\par Put $f(z)=(Pz, g(z))$ for $z\in V$. By Theorem 2.8 $f$ is ${\cal
A}_p^n$-biholomorphic in some neighborhood $U\subset V$ of $\zeta $.
Then $X\cap U = \{ z\in U: f_{k+1}(z)=0,...,f_n(z)=0 \} $, since
$(f_{k+1},...,f_n)=g$ and $X\cap U = \{ z\in U: g(z)=0 \} $.
\par {\bf 2.10. Definitions.} Let $U$ be an open subset in ${\cal A}_p^n$,
$2\le p<\infty $ or $p=\Lambda $. A subset $X$ in $U$ is called a
${\cal A}_p$-submanifold of ${\cal A}_p^n$ if the equivalent
conditions of Corollary 2.9 are satisfied. If in addition $X$ is a
closed subset in $U$, then $X$ is called a closed ${\cal
A}_p$-submanifold of $U$. This definition is the particular case of
the following general definition.
\par An ${\cal A}_p$-holomorphic manifold of ${\cal A}_p$-dimension $n$ is
a real $2^pn$-dimensional or $card (\Lambda )$-dimensional
$C^{\infty }$-manifold $X$ together with a family $ \{ (U_j,\phi
_j): j\in \Psi \} $ of charts such that
\par $(i)$ each $U_j$ is an open subset in $X$ and
$\bigcup_{j\in \Psi }U_j=X$, where $\Psi $ is a set;
\par $(ii)$ for each $j\in \Psi $ a mapping
$\phi _j: U_j\to V_j$ is a homeomorphism on an open subset $V_j$ in
${\cal A}_p^n$;
\par $(iii)$ for each $j, l\in \Psi $ a connection mapping
$\phi _j\circ \phi _l^{-1}$ is an ${\cal A}_p$-biholomorphic map
(see \S 2.6) from $\phi _l(U_j\cap U_l)$ onto $\phi _j(U_j\cap
U_l)$. Such system is called an ${\cal A}_p$-holomorphic atlas
$At(X):= \{ (U_j,\phi _j): j\in \Psi \} $. Each chart $(U_j,\phi
_j)$ provides a system of ${\cal A}_p$-holomorphic coordinates
induced from ${\cal A}_p^n$. For short we shall write ${\cal
A}_p$-manifold instead of ${\cal A}_p$-holomorphic manifold and
${\cal A}_p$-atlas instead of ${\cal A}_p$-holomorphic atlas if
other will not be specified.
\par For two ${\cal A}_p$-manifolds $X$ and $Y$
with atlases $At(X):= \{ (U_j,\phi _j): j\in \Psi _X \} $ and
$At(Y):= \{ (W_l,\psi _l): l\in \Psi _Y \} $ a function $f: X\to Y$
is called ${\cal A}_p$-holomorphic if $\psi _l\circ f\circ \phi
_j^{-1}$ is ${\cal A}_p$-holomorphic on $\phi _j(U_j\cap
f^{-1}(W_l))$. If $f: X\to Y$ is an ${\cal A}_p$-biholomorphic
epimorphism, then $X$ and $Y$ are called ${\cal
A}_p$-biholomorphically equivalent.
\par A subset $Z$ of an ${\cal A}_p$-manifold $X$ is called
an ${\cal A}_p$-submanifold, if $\phi _j(U_j\cap Z)$ is an ${\cal
A}_p$-submanifold in ${\cal A}_p^n$ for each chart $(U_j,\phi _j)$.
If additionally $Z$ is closed in $X$, then $Z$ is called a closed
${\cal A}_p$-submanifold.
\par {\bf 2.11. Theorem.} {\it Let $n\ge 2$, $f_1,...,f_n
\in C^1_{0,(z,\tilde z)}({\cal A}_p^n,{\cal A}_p)$ with $2\le p
<\infty $ or $p=\Lambda $ be a family of continuously ${\cal A}_p$
$(z,\tilde z)$-superdifferentiable functions satisfying
compatibility conditions:
$$(i)\quad \partial f_j/\partial {\tilde z}_k=
\partial f_k/\partial {\tilde z}_j\mbox{ for each }
j, k=1,...,n,$$ where $C^1_{0,(z,\tilde z)}({\cal A}_p^n,{\cal
A}_p)$ is the subspace of $C^1_{(z,\tilde z)}({\cal A}_p^n,{\cal
A}_p)$ of functions with closed bounded support. Then there exists
$u\in C^1_{0,(z,\tilde z)}({\cal A}_p^n,{\cal A}_p)$ satisfying the
following $\tilde \partial $-equation:
$$(ii)\quad \partial u/\partial {\tilde z}_j={\hat f}_j,
\quad j=1,...,n; $$ in particular, $(\partial u/\partial {\tilde
z}_j).1=f_j.$}
\par {\bf Proof.} Using the beginning and the end of the proof of
Theorem 2.3 we reduce the proof of Theorem 2.11 to the case of
finite $p$ mentioning, that the intersection $A\cap \theta ({\cal
A}_p)^n$ of a closed bounded subset $A$ in ${\cal A}_{\Lambda }^n$
for finite $p$ is compact. We put
$$(iii)\quad u(z) := - (2\pi )^{1-2^p}\int_{\theta ({\cal A}_p)}
[({\hat f}_1(\zeta _1,z_2,...,z_n).d{\tilde \zeta }_{2^p}) \wedge
\eta ],\mbox{ where}$$
$$\eta :=(\partial _{\zeta _1} Ln(\zeta _1-\zeta _2))M_1^*]). [(\partial
_{\zeta _2}Ln (\zeta _2-\zeta _3))M_2^*])...). [(\partial _{\zeta
_m}Ln (\zeta _m-z))M_m^*]),$$ $m=2^p-1$ (see \S 2.3). By changing of
variables we get
$$u(z):=-(2\pi )^{1-2^p}\int_{{\cal A}_p}[({\hat f}_1(z_1+\zeta _1,z_2,...,z_n).
d{\tilde \zeta }_{2^p})\wedge \psi ],\mbox{ where}$$
$$\psi :=(\partial _{\zeta _1} Ln(\zeta _1-\zeta _2))M_1^*]). [(\partial
_{\zeta _2}Ln (\zeta _2-\zeta _3))M_2^*])...). [(\partial _{\zeta
_m}Ln (\zeta _m))M_m^*]),$$  $m=2^p-1$. Therefore, $u\in
C^1_{(z,\tilde z)}({\cal A}_p^n,{\cal A}_p)$. Due to Theorem 2.3
$\partial u/\partial {\tilde z}_1={\hat f}_1$ in ${\cal A}_p^n$. In
view of Theorem 2.1 and the condition $\partial f_1/\partial {\tilde
z}_k=\partial f_k/\partial {\tilde z}_1$ the following equality is
satisfied
$${\hat f}_k(z)=-(2\pi )^{1-2^p}\int_{{\cal A}_p} \{
[\partial {\hat f}_k(\zeta _1,z_2,...,z_n)/\partial {\tilde \zeta
}_1]. d{\tilde \zeta }_{2^p} \} \wedge \psi ,$$ hence $\partial
u/\partial {\tilde z}_k={\hat f}_k$ for $k=2,...,n$, that is, $u$
satisfies equations $(ii)$. From this it follows, that $u$ is ${\cal
A}_p$-holomorphic in ${\cal A}_p^n\setminus (supp (f_1)\cup ... \cup
supp (f_n))$. In view of formula $(iii)$ it follows, that there
exists $0<r<\infty $ such that
\par $(iv) \quad u(z)=0$ for each $z\in {\cal A}_p^n$
with $|z_2|+...+|z_n|>r$. From $\partial u/\partial {\tilde
z}_1={\hat f}_1$ it follows, that $\partial u/\partial {\tilde
z}_1=0$ in ${\cal A}_p^n\setminus supp (f_1)$. Consequently, there
exists $0<R<\infty $ such that $u$ may differ from $0$ on ${\cal
A}_p^n\setminus B({\cal A}_p^n,0,R)$ only on an ${\cal A}_p$
constant (see Theorem 3.28 and Note 3.11 in \cite{luoystoc}).
Together with $(iv)$ this gives, that $u(z)=0$ on ${\cal
A}_p^n\setminus B({\cal A}_p^n,0,\max (R,r)).$
\par {\bf 2.12. Theorem.} {\it Let $U$ be an open subset in ${\cal A}_p^n$,
where $n\ge 2$, $2\le p<\infty $ or $p=\Lambda $. Suppose $K$ is a
bounded closed subset in $U$ such that $U\setminus K$ is connected.
Then for every ${\cal A}_p$-holomorphic function $h$ on $U\setminus
K$ there exists a function $H$ ${\cal A}_p$-holomorphic in $U$ such
that $H=h$ in $U\setminus K$.}
\par {\bf Proof.} Take any infinite $(z,\tilde z)$-differentiable
function $\chi $ on $U$ with bounded closed support such that $\chi
|_V=1$ on some (open) neighbourhood $V$ of $K$. Then consider a
family of functions $f_j$ such that ${\hat f}_j(z,\tilde z).S= -\{
(\partial \chi /\partial {\tilde z}).S \} h$ in $U\setminus K$ and
$f_j=0$ outside $U\setminus K$ for each $S$ in the set of generators
of ${\cal A}_p$, where $j=1,...,n$, $f_j(z)={\hat f}_j(z).1$.
Therefore, conditions of Theorem 2.11 are satisfied and it gives a
function $u\in C^1_{0,(z,\tilde z)}({\cal A}_p^n,{\cal A}_p)$ such
that $\partial u/\partial {\tilde z}_j={\hat f}_j$ for each
$j=1,...,n$. A desired function $H$ can be defined by the formula
$H:=(1-\chi )h-u$ such that $H$ is ${\cal A}_p$-holomorphic in $U$.
Since $\chi $ has a bounded closed support, then there exists an
unbounded connected subset $W$ in ${\cal A}_p^n\setminus supp (\chi
)$. Therefore, $u|_W=0$, consequently, $H|_{U\cap W}=h|_{U\cap W}$.
From $(U\setminus K)\cap W\ne \emptyset $ and connectedness of
$U\setminus K$ it follows, that $H|_{U\setminus K}=h|_{U\setminus
K}$.
\par {\bf 2.13. Remark.} In the particular case of a singleton
$K=\{ z \} $ Theorem 2.12 gives nonexistence of isolated
singularities, that is, each ${\cal A}_p$-holomorphic function in
$U\setminus \{ z \} $ for $U$ open in ${\cal A}_p^n$ with $n\ge 2$
can be ${\cal A}_p$-holomorphically extended to $z$. Theorem 2.12 is
the ${\cal A}_p$-analog of the Hartog's theorem for $\bf C^n$.
\par {\bf 2.14. Corollary.} {\it Let $U$ be an open
connected subset in ${\cal A}_p^n$, $2\le p<\infty $ or $p=\Lambda $
and $n\ge 2$. Suppose that $f$ is a right superlinearly ${\cal
A}_p$-superdifferentiable function $f: U\to {\cal A}_p$ and
$N(f):=\{ z\in U: f(z)=0 \} $, then
\par $(i)$ $U\setminus N(f)$ is connected,
\par $(ii)$ $N(f)$ is not bounded closed.}
\par {\bf Proof.} Reduce the proof of this corollary to the case of finite $p$ using
the beginning and the end of the proof of Theorem 2.3. Since if
statements $(i,ii)$ are true in the projection from ${\cal
A}_{\Lambda }$ on $\theta ({\cal A}_p)$, then they are true for
${\cal A}_{\Lambda }$. $(i)$. Write $f$ in the form
$f=\sum_{l=1}^{2^{p-1}}i_lg_l$, where
$g_l:=f_{2l-1}+i_{2l-1}^*i_{2l}f_{2l}$,
$f=\sum_{s=1}^{2^p}i_{s-1}f_s$, $f_s$ are real-valued functions, $\{
i_0,...,i_{2^p-1} \} $ is the set of generators of the
Cayley-Dickson algebra ${\cal A}_p$. In view of Proposition 2.3 and
Corollary 2.5.1 \cite{luoystoc} each function $g_l$ is holomorphic
in complex variables $y_k$, $k=1,...,2^{p-1}$, where
$z=\sum_{l=1}^{2^{p-1}}i_{2l-1}y_l$, $y_l=x_{2l-1}
+i_{2l-1}^*i_{2l}x_{2l}$, $x_1,...,x_{2^p}\in \bf R$, $z\in {\cal
A}_p$. \par Therefore, $N(f)=\bigcap_{l=1}^{2^{p-1}}N(g_l)$,
consequently, $U\setminus N(f)=\bigcup_{l=1}^{2^{p-1}}(U\setminus
N(g_l)).$ Then from Corollary 1.2.4 \cite{henlei} for complex
holomorphic functions $(i)$ follows.
\par $(ii)$. Suppose that $N(f)$ is bounded closed (compact for finite $p$).
In view of $(i)$ and Theorem 2.12 the function $1/f$ can be ${\cal
A}_p$-holomorphically extended on $N(f)$. This is the contradiction,
since $f=0$ on $N(f)$.
\par {\bf 2.14.1. Note.} Corollary $2.14$ is not true for arbitrary
${\cal A}_p$-holomorphic functions, for example, $f(\mbox{
}^1z,\mbox{ }^2z)=f_1(\mbox{ }^1z)f_2(\mbox{ }^2z)$ on $B({\cal
A}_p^2,0,2)$, where $f_1(\mbox{ }^1z):= - \mbox{
}^1z(\sum_{l=1}^{2^p}i_{l-1}\mbox{ }^1zi_{l-1})/(2^p-2) -r_1$,
$f_2(\mbox{ }^2z):= - \mbox{ }^2z(\sum_{l=1}^{2^p}i_{l-1}\mbox{ }^2
zi_{l-1})/(2^p-2) - r_2$ for finite $2\le p$, $0<r_1$, $0<r_2$,
$r_1^2+r_2^2<4$.
\par {\bf 2.15. Theorem.} {\it Let $U$ be an open subset
in ${\cal A}_p^n$, $2\le p<\infty $ or $p=\Lambda $, $f_1$,...,$f_n$
be infinite Frech\'et differentiable (by real variables) functions
on $U$ and suppose $(z,{\tilde z})$-superdifferentiable that
Conditions $2.11.(i)$ are satisfied in $U$. Then for each open
bounded polytor $P=P_1\times ...\times P_n$ such that $cl (P)$ is a
subset in $U$, there exists a function $u$ infinite differentiable
(by real variables) on $P$ and satisfying Conditions $2.11.(ii)$ on
$P$.}
\par {\bf Proof.} Using the beginning and the end of the proof of Theorem 2.3
reduce the proof of this theorem to finite $p$. Suppose that the
theorem is true for $f_{m+1}=...=f_n=0$ on $U$. The case $m=0$ is
trivial. Assume that the theorem is proved for $m-1$. Consider
$U'={U'}_1\times ...\times {U'}_n$ and $U"={U"}_1\times ...\times
{U"}_n$ open polytors in ${\cal A}_p^n$ such that $P\subset cl
(P)\subset {U"}\subset cl ({U"}) \subset {U'}\subset cl (U')\subset
U$. Take an infinite differentiable (by real variables) function
$\chi $ on ${U'}_m$ with compact support such that $\chi
|_{{U"}_m}=1$, $\chi =0$ in a neighbourhood of ${\cal A}_p\setminus
{U'}_m$. There exists a function
$$\eta (z):=-(2\pi )^{1-2^p}\int_{{U'}_m} [ \chi (\zeta )
({\hat f}_m(\mbox{ }^1z,...,\mbox{ }^{m-1}z,\zeta _1, \mbox{
}^{m+1}z,...,\mbox{ }^nz).d{\tilde \zeta }_{2^p})] \wedge \nu ,$$
where a differential form $\nu $ is given in \S 2.3 with $\zeta _1,
\zeta _2,..., \zeta _{2^p-1} \in {U'}_m$ and $\mbox{ }^mz$ here for
$\nu $ instead of $z$ in \S 2.3. By changing of variables as in \S
2.3 we get
$$\eta (z):=-(2\pi )^{1-2^p}\int_{{\cal A}_p}[\chi (\zeta _1+z)
({\hat f}_m(\mbox{ }^1z,...,\mbox{ }^{m-1}z, \zeta _1+\mbox{
}^mz,\mbox{ }^{m+1}z,...,\mbox{ }^nz ).d{\tilde \zeta
}_{2^p})]\wedge \psi ,$$ where the differential form $\psi $ is the
same as in \S 2.11. Consequently, $\partial \eta /\partial {\tilde
z}_m={\hat f}_m$ in ${U"}$. In view of Conditions 2.11(i) and
differentiating under the sign of the integral, since the support of
$\chi $ is compact, we get $\partial \eta (z)/\partial {\tilde
z}={\hat f}_j=0$ on $U'$ for $j=m+1,...,n$, since $f_j=0$, ${\hat
f}_j.1=f_j$, ${\hat f}_j$ is the partial (super)derivative of some
function $\xi _j$ by $z$. Thus functions $g_j:=f_j-\partial \eta
/\partial {\tilde z}$ for $j=1,...,n$ fulfil the compatibility
conditions 2.11(i), consequently, $g_m,...,g_n=0$ in $U"$. And
inevitably by the induction hypothesis there exists a function $v\in
C^{\infty }(P,{\cal A}_p)$ such that $\partial v/\partial {\tilde
z}_j=g_j$ in $P$ for which $u=v+\eta $ is the required solution.
\par {\bf 2.16. Definition.} Let $W$ be an open subset in ${\cal
A}_p^n$, $2\le p<\infty $ or $p=\Lambda $ and for each open subsets
$U$ and $V$ in ${\cal A}_p^n$ such that
\par $(i)$ $\emptyset \ne U\subset V\cap W\ne V$ and
\par $(ii)$ $V$ is connected  \\
there exists an ${\cal A}_p$-holomorphic (right superlinearly
superdifferentiable, in short RSS, correspondingly) function $f$ in
$W$ such that there does not exist any ${\cal A}_p$-holomorphic
(RSS) function $g$ in $V$ such that $g=f$ in $U$. Then $W$ is called
a domain of ${\cal A}_p$ (RSS, respectively) holomorphy. Sets of
${\cal A}_p$-holomorphic (RSS) functions in $W$ are denoted by
${\cal H}(W)$ (${\cal H}_{RSS}(W)$ respectively).
\par {\bf 2.17. Definition.} Suppose that $W$ is an open subset
in ${\cal A}_p^n$, $2\le p<\infty $ or $p=\Lambda $ and $K$ is a
closed bounded subset of $W$, then
\par $(i)$ ${\hat K}^{\cal H}_W:= \{ z\in W: |f(z)| \le
\sup_{\zeta \in K}\| {\hat f}(\zeta ) \| $ $\mbox{for each}$ $f\in
{\cal H}(W) \} $;
\par $(ii)$ ${\hat K}^{{\cal H}_{RSS}}_W:= \{ z\in W: |f(z)| \le
\sup_{\zeta \in K}|f(\zeta )|$ $\mbox{for each}$
$f\in {\cal H}_{RSS}(W) \} $;  \\
these sets are called the ${\cal H}(W)$-convex hull of $K$ and the
${\cal H}_{RSS}(W)$-convex hull of $K$ respectively, where $\| {\hat
f}(\zeta )\| :=\sup_{h\in {{\cal A}_p^n}, |h|\le 1} |{\hat f}(\zeta
).h|$. If $K={\hat K}^{\cal H}_W$ or $K={\hat K}^{{\cal
H}_{RSS}}_W$, then $K$ is called ${\cal H}(W)$-convex or ${\cal
H}_{RSS}(W)$-convex correspondingly.
\par {\bf 2.18. Proposition.} {\it For each closed bounded subset $K$
in ${\cal A}_p^n$, $2\le p<\infty $ or $p=\Lambda $, the ${\cal
H}({\cal A}_p^n)$-hull and ${\cal H}_{RSS}( {\cal A}_p^n)$-hull of
$K$ are contained in the $\bf R$-convex hull of $K$.}
\par {\bf Proof.} Reduce the proof to the case of finite $p$ using \S 2.3.
\par $I.$ Consider at first the ${\cal H}({\cal A}_p^n)$-hull of $K$.
Each $z\in {\cal A}_p^n$ can be written in the form $z=(\mbox{
}^1z,...,\mbox{ }^nz)$, $\mbox{ }^jz\in {\cal A}_p$, $\mbox{
}^jz=\sum_{l=1}^{2^p}x_{l,j}S_l$, where $x_{l,j}=x_{l,j}(z)\in \bf
R$, $S_l=i_{l-1}$. If $w\in {\cal A}_p^n$, $w\notin co_{\bf R}(K)$,
then there are $y_1,...,y_{2^pn}\in \bf R$ such that
$\sum_{j=1}^n\sum_{l=1}^{2^p} x_{l,j}(w)y_{2^p(j-1)+l}=0$, but \\
$\sum_{j=1}^n\sum_{l=1}^{2^p}x_{l,j}(w) y_{2^p(j-1)+l}<0$ if $z\in K$, where \\
$co_{\bf R}(K):=\{ z\in {\cal A}_p^n:$ $\mbox{there are}$
$a_1,...,a_s\in \bf R$ $\mbox{and}$ $v_1,...,v_s\in K$ $\mbox{such
that}$ $z=a_1v_1+...+a_kv_k \} $ denotes a $\bf R$-convex hull of
$K$ in ${\cal A}_p^n$. Put $\zeta _j=\sum_{l,j}y_{2^p(j-1)+l}S_l$,
then $f(z):=\exp (\sum_{j=1}^nz_j{\tilde \zeta }_j)$ is the ${\cal
A}_p$-holomorphic function in ${\cal A}_p^n$ such that $|f(z)|<1$
for each $z\in K$ and $|f(w)|=1$ for the marked point $w$ above (see
Corollary 3.3 \cite{luoyst}), since $i_v^2 = -1$ for each $v>0$.
From $\| {\hat f}(\zeta ) \| \ge |f(\zeta )|$ the first statement
follows.
\par $II.$ Consider now the ${\cal H}_{RSS}({\cal A}_p^n)$-hull
of $K$. Each $f\in {\cal H}_{RSS}(W)$ has the form
$f=\sum_{l=1}^{2^{p-1}}i_lg_l$, where each function $g_l$ is
holomorphic in complex variables $y_k$ (see \S 2.14). \par The set
$K$ has projection $K_k$ on the complex subspaces $\bf C^n$
corresponding to variables $\mbox{ }^1y_k,...,\mbox{ }^ny_k$.
Therefore, $({\hat K}^{{\cal H}_{RSS}}_{{\cal A}_p^n})_k \subset
{{\hat K}_{k,\bf C^n}}^{\cal O}$ for each $k$, where ${{\hat
K}_{k,\bf C^n}}^{\cal O}$ denotes the complex holomorphic hull of
$K_k$ in $\bf C^n$. In view of Proposition 1.3.3 \cite{henlei}
${{\hat K}_{k,\bf C^n}}^{\cal O}\subset co_{\bf R}(K_k)$, hence
${\hat K}^{{\cal H}_{RSS}}_{{\cal A}_p^n}\subset co_{\bf R}(K)$.
\par {\bf 2.18.1. Note.} Due to Proposition $2.18$ above Corollary $1.3.4$ \cite{henlei}
can be transferred on $\cal H$ and ${\cal H}_{RSS}$ for ${\cal
A}_p^n$ instead of $\bf C^n$. Also ${\cal A}_p$-versions of Theorems
$1.3.5, 7, 11$, Corollaries $1.3.6, 8, 9, 10, 13$ and Definition
$1.3.12$ are true in the ${\cal H}_{RSS}$-class of functions instead
of complex holomorphic functions.
\section{Integral representations of functions of Cayley-Dickson
variables}
\par {\bf 3.1. Definitions and Notations.} Consider an
${\cal A}_p$-valued function on ${\cal A}_p^n$, $2\le p<\infty $
or $p=\Lambda $ such that \\
$(i)\quad (\zeta ,\zeta )=ae$ with
$a\ge 0$ and $(\zeta ,\zeta )=0$ if and only if $\zeta =0$, \\
$(ii)\quad (\zeta ,z+\xi )=(\zeta ,z)+(\zeta ,\xi )$, \\
$(iii)\quad (\zeta +\xi ,z)=(\zeta ,z)+(\xi ,z)$, \\
$(iv)\quad (\alpha \zeta ,z)=\alpha (\zeta ,z)=(\zeta ,\alpha z)$
for each $\alpha \in \bf R$ and $(\zeta \alpha ,\zeta )={\tilde
{\alpha }}(\zeta ,\zeta )$ for each $\alpha \in {\cal A}_p$, \\
$(v)\quad (\zeta ,z)^{\tilde .}=(z,\zeta )$ for each $\zeta , \xi $
and $z\in {\cal A}_p^n$, $n\in \bf N$. Then this function is called
the scalar product in ${\cal A}_p^n$.
The corresponding norm is: \\
$(vi)\quad |\zeta |=\{ (\zeta ,\zeta ) \} ^{1/2}$.
In particular, it is possible to take the canonical scalar product: \\
$(vii)\quad <\zeta ;z>:= (\zeta ,z)=\sum_{l=1}^n\mbox{ }^l{\tilde
\zeta }\mbox{ }^lz$, where $z=(\mbox{ }^1z,...,\mbox{ }^nz)$,
$\mbox{ }^lz\in {\cal A}_p$.
\par Consider differential forms on ${\cal A}_p$: \\
$(1)\quad \phi _{p,0}(z) := d{\tilde z}\wedge d{\tilde z},$ ${\phi
'}_{p,0}(z):={\tilde z}d{\tilde z}$, \\
$\phi _{p,k}(z) :=(i_{2k}(d{\tilde z}i_{2k}))\wedge (i_{2k}(dz
i_{2k})),$ ${\phi '}_{p,k}(z) :=(i_{2k}({\tilde z}i_{2k}))\wedge
(i_{2k}(dz i_{2k})),$ for each $k=1,...,2^{p-1}-1$, \\
$(2)\quad w_{2^p}(z) := C_p \{ \phi _{p,0}(z)\wedge \phi
_{p,1}(z)\wedge ... \wedge \phi _{p,2^{p-1}-1}(z) \} _{q_0(2^{p-1})}$, \\
where $C_p=const \ne 0$; \\
$(3)\quad w_{2^p,k}(\zeta -z):= \{ \phi _{p,0}(\zeta )\wedge ...
\wedge \phi _{p,k-1}(\zeta )\wedge {\phi '}_{p,k} (\zeta - z)\wedge
\phi _{p,k+1}(\zeta )\wedge \phi _{p,2^{p-1}-1}(\zeta ) \}
_{q_0(2^{p-1})}$
for each $k=0,...,2^{p-1}-1$, \\
where $q_0(s)$ means the associated product in definite order
corresponding to the left preferred order of brackets, $\{ b_1...b_s
\} _{q_0(s)}:=(...((b_1b_2)b_3)...b_{s-1})b_s$ for $b_1,...,b_s\in
{\cal A}_p$. Introduce also differential forms: \\
$(4)\quad {\check{\phi }}_{p,0}(\zeta ,z) := (d{\tilde {\zeta
}}-d{\tilde z})\wedge d{\tilde {\zeta }},$  ${\check{\phi
}'}_{p,0}(\zeta ,z):=({\tilde {\zeta }}- {\tilde z}) d{\tilde {\zeta }}$, \\
${\check{\phi }}_{p,k}(\zeta ,z) :=(i_{2k}(d{\tilde {\zeta
}}-d{\tilde z})i_{2k}))\wedge (i_{2k}(d\zeta  i_{2k})),$
${\check{\phi }'}_{p,k}(\zeta ,z) :=(i_{2k}(({\tilde {\zeta
}}-{\tilde z})i_{2k}))\wedge
(i_{2k}(d\zeta i_{2k})),$ for each $k=1,...,2^{p-1}-1$, \\
$(5)\quad {\check{w}}_{2^p}(\zeta ,z) := C_p \{ {\check{\phi
}}_{p,0}(\zeta ,z)\wedge {\check{\phi }}_{p,1}(\zeta ,z)\wedge ...
\wedge {\check{\phi
}}_{p,2^{p-1}-1}(\zeta ,z) \} _{q_0(2^{p-1})}$, \\
where $C_p=([(2^p-2)!]2(2^{p-1}-1))^{-1}$; \\
$(6)\quad {\check{w}}_{2^p,k}(\zeta ,z):= \{ {\check{\phi
}}_{p,0}(\zeta ,z)\wedge ... \wedge {\check{\phi }}_{p,k-1}(\zeta
,z)\wedge {\check{\phi }'}_{p,k} (\zeta, z)\wedge {\check{\phi
}}_{p,k+1}(\zeta ,z)\wedge ... \wedge {\check{\phi
}}_{p,2^{p-1}-1}(\zeta ,z) \} _{q_0(2^{p-1})}$
for each $k=0,...,2^{p-1}-1$, \\
$(7)\quad {\hat {\phi }}_{p,0}(\zeta ,z) := (d{\tilde {\zeta
}}-d{\tilde z})\wedge (d{\tilde {\zeta }}-d{\tilde z}),$ ${\hat
{\phi }'}_{p,0}(\zeta ,z):=({\tilde {\zeta }}- {\tilde z})
(d{\tilde {\zeta }} - d{\tilde z})$, \\
${\hat {\phi }}_{p,k} (\zeta ,z) :=(i_{2k}(d{\tilde {\zeta
}}-d{\tilde z})i_{2k}))\wedge (i_{2k}((d\zeta  -dz)i_{2k}))),$
${\hat{\phi }'}_{p,k}(\zeta ,z) :=(i_{2k}(({\tilde {\zeta }}-{\tilde
z})i_{2k}))\wedge (i_{2k}((d\zeta -dz)i_{2k})),$ for each $k=1,...,2^{p-1}-1$, \\
$(8)\quad {\hat w}_{2^p}(\zeta ,z) := C_p \{ {\hat {\phi
}}_{p,0}(\zeta ,z)\wedge {\hat {\phi }}_{p,1}(\zeta ,z)\wedge ...
\wedge {\hat {\phi }}_{p,2^{p-1}-1}(\zeta ,z) \} _{q_0(2^{p-1})}$, \\
where $C_p=([(2^p-2)!]2(2^{p-1}-1))^{-1}$; \\
$(9)\quad {\hat w}_{2^p,k}(\zeta ,z):= \{ {\hat {\phi }}_{p,0}(\zeta
,z)\wedge ... \wedge {\hat {\phi }}_{p,k-1}(\zeta ,z)\wedge {\hat
{\phi }'}_{p,k} (\zeta, z)\wedge {\hat {\phi }}_{p,k+1}(\zeta
,z)\wedge ... \wedge {\hat {\phi }}_{p,2^{p-1}-1}(\zeta ,z) \}
_{q_0(2^{p-1})}$
for each $k=0,...,2^{p-1}-1$, \\
where we can express $\tilde {\zeta }$ and $\tilde z$ in the $\zeta
$ and $z$-representations respectively: ${\tilde z}=(2^p-2)^{-1} \{
-z + \sum_{s\in {\hat b}} s(z {\tilde s}) \} $ for each $2\le p\in
\bf N$. With the help of them construct differential forms on ${\cal A}_p^n$: \\
$$(10)\quad \theta _z(\zeta ):={C'}_p
|\zeta -z|^{-2^pn} \sum_{s=1}^n \sum_{q=0}^{2^{p-1}-1}\{
w_{2^p}(\mbox{ }^1\zeta )\wedge ...$$
$$\wedge w_{2^p}(\mbox{ }^{s-1}\zeta )\wedge w_{2^p,q}
(\mbox{ }^s\zeta -\mbox{ }^sz)\wedge w_{2^p}(\mbox{ }^{s+1}\zeta )
\wedge ... \wedge w_{2^p}(\mbox{ }^n\zeta ) \} _{q_0(n)},$$
$$(11)\quad {\check{\theta }}(\zeta ,z):={C'}_p |\zeta -z|^{-2^pn}
\sum_{s=1}^n \sum_{q=0}^{2^{p-1}-1} \{ {\check{w}}_{2^p}(\mbox{
}^1\zeta ,\mbox{ }^1z) \wedge ...$$
$$\wedge {\check{w}}_{2^p}(\mbox{ }^{s-1}\zeta ,\mbox{ }^{s-1}z)\wedge
{\check{w}}_{2^p,q}(\mbox{ }^s\zeta ,\mbox{ }^sz)\wedge
{\check{w}}_{2^p}(\mbox{ }^{s+1}\zeta ,\mbox{ }^{s+1}z)\wedge ...
\wedge {\check{w}}_{2^p}(\mbox{ }^n\zeta ,\mbox{ }^nz) \}_{q_0(n)}
;$$
$$(12)\quad {\hat {\theta }}(\zeta ,z):={C'}_p |\zeta -z|^{-2^pn}
\sum_{s=1}^n \sum_{q=0}^{2^{p-1}-1} \{ {\hat w}_{2^p}(\mbox{
}^1\zeta ,\mbox{ }^1z) \wedge ...$$
$$\wedge {\hat w}_{2^p}(\mbox{ }^{s-1}\zeta ,\mbox{ }^{s-1}z)\wedge
{\hat w}_{2^p,q}(\mbox{ }^s\zeta ,\mbox{ }^sz)\wedge {\hat
w}_{2^p}(\mbox{ }^{s+1}\zeta ,\mbox{ }^{s+1}z)\wedge ... \wedge
{\hat w}_{2^p}(\mbox{ }^n\zeta ,\mbox{ }^nz) \}_{q_0(n)} ,$$ where
${C'}_p:=(2^pn)!! (2\pi )^{-2^{p-1}n}$; $|\zeta -z|^2$ is considered
in the $(\zeta -z, {\tilde \zeta } -{\tilde z})$-representation:
$|\mbox{ }^s\zeta - \mbox{ }^sz|^2= (\mbox{ }^s\zeta -\mbox{
}^sz)(\mbox{ }^s{\tilde {\zeta }} - \mbox{ }^s{\tilde z})$, $|\zeta
-z|^2=\sum_{s=1}^n|\mbox{ }^s\zeta -\mbox{ }^sz|^2$, $\zeta $ and
$z\in {\cal A}_p^n$. If $U$ is an open subset in ${\cal A}_p^n$ and
$f$ is a bounded ${\cal A}_p$-differential form on $U$, then by the
definition:
$$(13)\quad ({\cal B}_Uf)(z):=\int_{\zeta \in U}f(\zeta )
\wedge {\check{\theta }}(\zeta ,z)$$ for each $z\in {\cal A}_p^n$.
If in addition $U$ is with a continuous piecewise $C^1$-boundary (by
the corresponding real variables) and $f$ is a bounded differential
form on $\partial U$, then by the definition:
$$(14)\quad ({\cal B}_{\partial U}f)(z):=
\int_{\zeta \in \partial U}f(\zeta )\wedge {\check{\theta }}(\zeta
,z)$$ for each $z\in {\cal A}_p^n$.
\par {\bf 3.2. Theorem.} {\it Let  $U$ be an open subset in ${\cal
A}_p^n$, $2\le p\in \bf N$, with piecewise $C^1$-boundary $\partial
U$. Suppose that $f$ is a continuous function on $cl (U)$ and
${\tilde \partial }f$ is continuous on $U$ in the sense of
distributions and has a continuous extension on $cl (U)$. Then
$$(1)\quad f={\cal B}_{\partial U}f-{\cal B}_U{\tilde \partial f}
\mbox{ on } U,$$ where ${\cal B}_U$ and ${\cal B}_{\partial U}$ are
the ${\cal A}_p$-integral operators given by Equations
$3.1.(13,14)$.}
\par {\bf Proof.} Write the variable $z$ in the form
$z = \sum_{l=0}^{2^{p-1}-1} i_{2l}\alpha _l$, where $\alpha _l\in
{\bf C}_l := {\bf R}\oplus i_{2l}^*i_{2l+1}\bf R$, $i_{2l}\alpha
_l = i_{2l}x_l + i_{2l+1}y_l$, where $x_l, y_l\in \bf R$. Then \\
$(1)\quad \alpha _0i_k=i_k{\bar {\alpha }}_0$ for each $k>0$, \\
$(2)\quad i_{2k}(dz i_{2k})=(\sum_{l>0, l\ne k} i_{2l}d\alpha _l)-
i_{2k}d{\bar {\alpha }}_k - d{\bar {\alpha }}_0$, \\
$(3)\quad (i_{2l}d\alpha _l)\wedge d{\bar {\alpha }}_0 = -d{\alpha }_0
\wedge (i_{2l}d\alpha _l)$ for each $l>0$, \\
$(4)\quad (i_{2l}d\alpha _l)\wedge (i_{2q}d{\alpha }_q) = (i_{2q}d
{\alpha }_q)\wedge (i_{2l}d\alpha _l)$ for each $l\ne q$ with $l>0$ and $q>0$, \\
$(5)\quad (i_{2l}d\alpha _l)\wedge (i_{2l}d{\bar {\alpha }}_l) =0$
for each $l>0$, then \\
$(6)\quad d\alpha _0\wedge d{\bar {\alpha }}_0 = -d{\bar {\alpha
}}_0\wedge d\alpha _0= -2i_1 dx_0\wedge dy_0$, $d\alpha _0\wedge d\alpha _0 =0$, \\
$(7)\quad (i_{2l}d\alpha _l)\wedge (i_{2l}d\alpha _l) =
2i_{2l}i_{2l+1} dx_l\wedge dy_l$, $d\alpha _l\wedge d{\bar {\alpha
}}_l = - 2i_{2l}^*i_{2l+1} dx_l\wedge dy_l=2i_{2l}i_{2l+1}
dx_l\wedge dy_l$ for each $l>0$.
\par From Equations $(1-7)$ and $d{\tilde z} = d{\bar {\alpha }}_0
- \sum_{l>0}i_{2l}d\alpha _l$ it follows, that \\
$(8)\quad \phi _{p,0}(z) = (d\alpha _0- d{\bar {\alpha }}_0)\wedge
(\sum_{q>0}i_{2q}d\alpha _q) + \sum_{l>0} \sum_{q>0}
(i_{2l}d\alpha _l)\wedge (i_{2q}d\alpha _q)$, \\
$(9)\quad \phi _{p,k}(z) = (d\alpha _0 \wedge d{\bar {\alpha }}_0) -
(d\alpha _k \wedge d{\bar {\alpha }}_k) + ((2i_{2k}d{\bar {\alpha
}}_k) \wedge (\sum_{l>0, l\ne k}i_{2l}d\alpha _l)) -  \sum_{l>0,
l\ne k} \sum_{q>0, q\ne k}(i_{2l}d\alpha _l)\wedge (i_{2q}d\alpha _q)$ \\
for each $k>0$. The differential form $w_{2^p}$ is of degree $2^p$
in real coordinates $x_0, y_0,...,$  $x_{2^{p-1}-1}, y_{2^{p-1}-1}$,
hence it may contain only the multiplier $dx_0\wedge dy_0$ or may
contain only $d\alpha _0 \wedge d{\bar {\alpha }}_0$, hence all
terms in $w_{2^p}$ arising from the term $(d\alpha _0 - d{\bar
{\alpha }}_0)\wedge (\sum_{q>0}i_{2q}d\alpha _q)$ in $\phi
_{p,0}(z)$ cancel, since $\phi _{p,k}(z)$ contains $d\alpha _0
\wedge d{\bar {\alpha }}_0$ for each $k>0$. Then from $(3-5)$ it
follows, that all terms arising from the term $((2i_{2k}d{\bar
{\alpha }}_k) \wedge (\sum_{l>0, l\ne k}i_{2l}d\alpha _l))$ in $\phi
_{p,k}(z)$ for $k>0$ cancel in $w_{2^p}$. Thus for a choice of the
multiplier $d\alpha _0 \wedge d{\bar {\alpha }}_0$ in $w_{2^p}$
there are $(2^{p-1}-1)$ possibilities among $\phi _{p,q}(z)$ with
$q=1,..,2^{p-1}-1$ in the graded external product. After a choice of
$d\alpha _0 \wedge d{\bar {\alpha }}_0$ for some $q>0$ it remains
$(2^p-2)(2^p-3)/2$ variants for a choice of the multiplier $d\alpha
_1 \wedge d{\bar {\alpha }}_1=(i_2d\alpha _1) \wedge (i_2d{\alpha
}_1)$. Then by induction after choices of the multipliers
$(i_{2v}d\alpha _v) \wedge (i_{2v}d{\alpha }_v)=d\alpha _v \wedge
d{\bar {\alpha }}_v$ for $v=0,1,...,q-1$ with $q>2$ it remains
$(2^p-2q)(2^p-2q-1)/2$ variants for choices of
the multiplier $(i_{2q}d\alpha _q) \wedge (i_{2q}d{\alpha }_q)$. Thus \\
$(10)\quad w_{2^p}=
(-1)^{2^{p-1}-2}C_p[(2^p-2)!](2^{p-1}-1)2^{-(2^{p-1}-1)} \{ (d\alpha
_0\wedge d{\bar {\alpha }}_0)\wedge ((i_2d\alpha _1)\wedge
(i_2d{\alpha }_1))\wedge ... \wedge ((i_{2^p-2}d\alpha
_{2^{p-1}-1})\wedge (i_{2^p-2}d{\alpha }_{2^{p-1}-1})) \}
_{q_0(2^{p-1})}$ \\
$= (C_p[(2^p-2)!]2(2^{p-1}-1)) dx_0\wedge dy_0\wedge dx_1\wedge
dy_1\wedge ... \wedge dx_{2^{p-1}-1}\wedge
dy_{2^{p-1}-1}$, \\
since $(i_0i_1)(i_2i_3)=i_1^2=-1$, $\{
(i_0i_1)(i_2i_3)...(i_{2^{p-2}}i_{2^{p-1}}) \} _{q_0(2^{p-1})}=-1$.
Hence $w_{2^p}$ is the volume element on ${\cal A}_p$ equal to the
Lebesgue measure $\mu $ on the underlying Euclidean space ${\bf
R}^{2^p}$ such that $\mu ([0,1]^{2^p})=1$, since $C_p=([(2^p-2)!]2
(2^{p-1}-1))^{-1}$.
\par The differential form $\check{\theta }(\zeta ,z)$ has the decomposition
$$(11)\quad {\check{\theta }}(\zeta ,z)=\sum_{q=0}^{2^{p-1}n-1}\Upsilon _q(\zeta ,z),$$
where $\Upsilon _q(\zeta ,z)$ is the ${\cal A}_p$-differential form
with all terms of degree $2^pn-q-1$ by $\zeta $ and $\tilde \zeta $
and their multiples on ${\cal A}_p$ constants and of degree $q$ by
$z$ and $\tilde z$ and their multiples on ${\cal A}_p$ constants.
The differential form $f(\zeta )$ has the decomposition
$$(12)\quad f(\zeta )=\sum_{r=0}^mf_r(\zeta ),$$ where
$m=deg (f)$, $f_r(\zeta )$ is with all terms of degree $r$ by $\zeta
$ and ${\tilde \zeta }$ and their multiples on ${\cal A}_p$
constants. Then $f_r\wedge \Upsilon _q=0$, when $r>q+1.$ By the
definition of integration $\int_{\zeta \in U}f_r(\zeta ) \wedge
\Upsilon _q(\zeta ,z)=0$ for $r<q+1.$ If $f$ is a function, then
$\int_{\zeta \in
\partial U} f(\zeta )\Upsilon _q(\zeta ,z)=0$ for each $q>0$, since
$\partial U$ has the dimension $2^pn-1$, hence
$$(13)\quad ({\cal B}_{\partial U}f)(z)=\int_{\zeta \in \partial U}
f(\zeta )\theta _z(\zeta ),$$ since $\Upsilon _0(\zeta ,z)=\theta
_z(\zeta )$. If $f$ is a $1$-form, then  $\int_{\zeta \in U}f(\zeta
) \wedge \Upsilon _q(\zeta ,z)=0$ for each $q>0$, since $U$ has the
dimension $2^pn$, consequently,
$$(14)\quad ({\cal B}_Uf)(z)=\int_{\zeta \in U} f(\zeta )\wedge \theta _z(\zeta ).$$
\par In particular, there are identities in $\bf H$:
$d\xi \wedge jd\zeta =d\xi j\wedge d\zeta $ and $(\xi d\zeta
)^{\tilde .}=[d{\tilde \zeta }]\tilde \xi $ for each
$\xi , \zeta \in \bf H$. Then  \\
$(i)\quad d\zeta \wedge jd{\tilde \zeta }\wedge
d{\tilde \zeta }\wedge d\zeta =0,$ \\
$(ii)\quad d\zeta \wedge d{\tilde \zeta }\wedge
jd{\tilde \zeta }\wedge d\zeta =0,$ \\
$(iii)\quad d\zeta \wedge d{\tilde \zeta }\wedge d{\tilde \zeta
}\wedge jd\zeta =0,$ since $j^2=-e$ and ${\bf R}e$ is the centre of
the quaternion algebra $\bf H$, $\alpha $ and $\beta \in \bf C$
commute with $d\alpha $, $d{\bar \alpha },$ $d\beta $ and $d{\bar
\beta }$, where $\zeta =\alpha e+\beta j$. From $(i-iii)$ with the
help of bijective surjective mappings $\zeta \mapsto
j\zeta $ and $\zeta \mapsto \zeta j$ it follows, that \\
$(iv)\quad d\zeta \wedge d{\tilde \zeta }\wedge
d{\tilde \zeta }\wedge d\zeta =0,$  \\
$(v)\quad d\zeta \wedge jd{\tilde \zeta }\wedge
jd{\tilde \zeta }\wedge d\zeta =0,$ \\
$(vi)\quad d\zeta \wedge jd{\tilde \zeta } \wedge jd{\tilde \zeta
}j\wedge jd\zeta j=0.$
\par Write $\xi \in {\cal A}_p$ in the form $\xi =\alpha  + \beta l$,
then ${\tilde \xi }={\tilde {\alpha }} - \beta l$, where $\alpha \in
{\cal A}_{p-1}$ and $\beta \in {\cal A}_{p-1}$, $l$ is the generator
of the doubling procedure of ${\cal A}_p$ from ${\cal A}_{p-1}$
\cite{baez}, since there is the identity $\beta l = l {\tilde {\beta
}}$. The decomposition $\xi =\alpha  + \beta l$ is unique for each
$\xi \in {\cal A}_p$, where $\alpha =\alpha (\xi )$ and $\beta
=\beta (\xi )$ depend on $\xi $ in general. Put \\
$\kappa _{p,0}(z) := dz\wedge d{\tilde z}$ and \\
$\kappa _{p,q}(z) := (i_{2q}(dz i_{2q}))\wedge (i_{2q}(dz i_{2q}))$
for each $1\le q\le 2^{p-1}-1$. From Formulas $(1-9)$ or from
Formulas $(i-vi)$ and induction
by $p$ with the help of doubling procedures it follows, that \\
$(15)\quad \{ \phi _{p,0}(z)\wedge ... \wedge \phi _{p,v-1}(z)\wedge
\kappa _{p,v}(z)\wedge \phi _{p,v+1}(z)\wedge ... \wedge \phi
_{p,2^{p-1}-1}(z) \}_{q_0(2^{p-1})}=0$ \\
for each $v = 0,1,...,2^{p-1}-1$.
\par In the $(\zeta -z, {\tilde {\zeta }} - {\tilde z})$-representation \\
$|\zeta -z|^2=\sum_{s=1}^n(\mbox{ }^s{\tilde {\zeta }} - \mbox{
}^s{\tilde z}) (\mbox{ }^s{\zeta } - \mbox{ }^sz)$, hence: \\
$(16)\quad d_{\zeta }|\zeta -z|^{2^pn} = (2^{p-1}n) |\zeta
-z|^{2^pn-2} \sum_{s=1}^n \{ (d\mbox{ }^s{\tilde {\zeta }}) (\mbox{
}^s{\zeta } - \mbox{ }^sz) + (\mbox{ }^s{\tilde {\zeta }} - \mbox{
}^s{\tilde z}) d\mbox{ }^s{\zeta } \} $. \\
From Formulas $(15,16)$ it follows, that \\
$(17)\quad d_{\zeta } (|\zeta -z|^{2^pn}\theta _z(\zeta ))= {C'}_p
2^{p-1}n \{ w_{2^p}(\mbox{ }^1{\zeta })\wedge ... \wedge
w_{2^p}(\mbox{ }^n{\zeta }) \} _{q_0(n)}$, \\
since $d_{\zeta }=\partial _{\zeta } +\partial _{\tilde {\zeta }}$.
Now calculate $d_{\zeta }\theta _z(\zeta )$ in $U\setminus \{ z \} $
using Formulas $(15-17)$: \\
$(18)\quad d_{\zeta }\theta _z(\zeta )=0$.
\par There exists $\epsilon _0>0$ such that for each
$0<\epsilon <\epsilon _0$ the ball $B({\cal A}_p^n,z,\epsilon ):= \{
\zeta \in {\cal A}_p^n: |\zeta -z|\le \epsilon \} $ and hence the
sphere $S({\cal A}_p^n,z,\epsilon ):= \{ \zeta \in {\cal A}_p^n:
|\zeta -z|=\epsilon \} =\partial B({\cal A}_p^n, z,\epsilon )$ are
contained in $U$. Apply the Stoke's formula for vector-valued
functions and differential forms componentwise, using the Euclidean
space $\bf R^{2^pn}$ underlying ${\cal A}_p^n$, then \\
$(19)\quad \int_{S({\cal A}_p^n,z,\epsilon )}f(\zeta )\theta
_z(\zeta )= \int_{\partial U}f(\zeta )\theta _z(\zeta
)-\int_{U_{\epsilon }} [df(\zeta )]\wedge \theta _z(\zeta )$, where
$U_{\epsilon }:=U\setminus B({\cal A}_p^n,z,\epsilon )$, $0<\epsilon
<\epsilon _0$. Therefore, from $(15,18)$ it follows, that
\par $(20)\quad {\cal B}_Udf={\cal B}_U{\tilde \partial }f$,
since $df=\partial f+{\tilde \partial f}$, where $\partial f(\zeta
)= (\partial f(\zeta )/\partial \zeta ).d\zeta $, ${\tilde \partial
}f(\zeta )=(\partial f(\zeta )/\partial {\tilde \zeta }). d{\tilde
\zeta }$, $f(\zeta )=f(\zeta ,{\tilde \zeta })$ is the abbreviated
notation.
\par In view of Formula $(18)$ and the Stoke's formula: \\
$(21)\quad \int_{S({\cal A}_p^n,z,\epsilon )}\theta _z(\zeta )= [
(2\pi )^{2^{p-1}n}/(2^pn)!!]^{-1} [\epsilon ^{-4n}]\int_{B({\cal
A}_p^n,z,\epsilon )} (dV)e=e$, where $dV$ is the standard volume
element of the Euclidean space $\bf R^{2^pn}$. In the even
dimensional Euclidean space ${\bf R}^{2k}$ the volume $V_{2k}$ of
the ball of radius $1$ relative to the standard Lebesgue measure
$\lambda $ with $\lambda ([0,1]^{2k})=1$ is $V_{2k}=(2\pi
)^k/(2k)!!$ (see \S XI.4.2, Example 3, in \cite{zorich}).
Then Formula $(21)$ implies, that \\
$\lim_{\epsilon \to 0}\int_{S({\cal A}_p^n,z,\epsilon )}
f(\zeta )\theta _z(\zeta )=f(z)$, since  \\
$\int_{S({\cal A}_p^n,z,\epsilon )}(f(\zeta )-f(z))\theta _z(\zeta
)= \epsilon ^{-2^pn+1}\int_{S({\cal A}_p^n,z,\epsilon )} (f(\zeta
)-f(z))[|\zeta -z|^{2^pn-1}\theta _z(\zeta )]$. The form $[|\zeta
-z|^{2^pn-1}\theta _z(\zeta )]$ is bounded on $U$, consequently,
$|\int_{S({\cal A}_p^n,z,\epsilon )} (f(\zeta )-f(z))\theta _z(\zeta
)|\le C_1\max \{ |f(\zeta )-f(z)|: \zeta \in B({\cal
A}_p^n,z,\epsilon ) \} $, where $C_1$ is a positive constant
independent of $f$ and $\epsilon $ for each $0<\epsilon <\epsilon
_0$. Therefore, Formula $(1)$ follows from Formula $(19)$ by taking
the limit when $\epsilon >0$ tends to zero and using Identity
$(20)$.
\par {\bf 3.3. Corollary.} {\it Let $U$ be an open subset in ${\cal
A}_p^n$, $2\le p\in \bf N$, and $f$ be a continuous function on $cl
(U)$ and ${\cal A}_p$-holomorphic on $U$. Then  \\
$(1)\quad f={\cal B}_{\partial U}f$ on $U$, \\
where ${\cal B}_U$ and ${\cal B}_{\partial U}$ are the integral
operators given by Equations $3.1.(13,14)$.}
\par {\bf Proof.} From ${\tilde \partial f}=0$, since
$\partial f(\zeta )/\partial {\tilde \zeta }=0$, and Formula
$3.2.(1)$ implies Formula $3.3.(1)$.
\par {\bf 3.4. Definitions and Notations.}
Suppose that $U$ is a bounded open subset in ${\cal A}_p^n$ and
$\psi (\zeta ,z)$ be an ${\cal A}_p$-valued $C^1$-function (by the
corresponding real variables) defined on $V\times U$,
where $V$ is a neighbourhood of $\partial U$ in ${\cal A}_p^n$, such that \\
$(1)\quad <\psi (\zeta ,z);\zeta -z>\ne 0$ for each $(\zeta ,z)\in
{\partial U}\times U$. Then $\psi $ is called an ${\cal
A}_p$-boundary distinguishing map.
Consider the function: \\
$(2)\quad \eta ^{\psi }(\zeta ,z,\lambda ):=
\lambda (\zeta -z) <\zeta -z;\zeta -z>^{-1}$ \\
$+(1-\lambda )\psi (\zeta ,z)<\zeta -z;\psi (\zeta ,z)>^{-1}$, \\
(see Formula $3.1.(vii)$) and the differential forms: \\
$(3)\quad {\check{\phi }}_{p,0}(\mbox{ }^s{\tilde \eta }^{\psi
}(\zeta ,z,\lambda ), \mbox{ }^s\zeta ) := [( {\tilde \partial
}_{\mbox{ }^s\zeta ,\mbox{ }^sz} + d_{\lambda })\mbox{ }^s{\tilde
\eta }^{\psi
} (\zeta ,z,\lambda )]\wedge d\mbox{ }^s {\tilde {\zeta }}$, \\
${\check{\phi }'}_{p,0}(\mbox{ }^s{\tilde \eta }^{\psi }(\zeta
,z,\lambda ), \mbox{ }^s\zeta ) := \mbox{ }^s{\tilde \eta }^{\psi
}(\zeta ,z,\lambda ) d \mbox{ }^s{\tilde {\zeta }}$, \\
$(4)\quad {\check{\phi }}_{p,u}(\mbox{ }^s{\tilde \eta }^{\psi
}(\zeta ,z,\lambda ), \mbox{ }^s\zeta ) := (i_{2u} \{ [({\tilde
\partial }_{\mbox{ }^s\zeta ,\mbox{ }^sz} +d_{\lambda })\mbox{
}^s{\tilde \eta }^{\psi }(\zeta ,z,\lambda )]i_{2u} \} ) \wedge
(i_{2u}(d \mbox{ }^s\zeta i_{2u}))$, \\
${\check{\phi }'}_{p,u}(\mbox{ }^s{\tilde \eta }^{\psi }(\zeta
,z,\lambda ), \mbox{ }^s\zeta ):= [i_{2u}(\mbox{ }^s{\tilde {\eta
}}^{\psi }(\zeta
,z,\lambda )i_{2u})] (i_{2u}(d\mbox{ }^s\zeta i_{2u}))$ for each $u>0$, \\
$(5)\quad {\check{\phi }}_{p,0}(\mbox{ }^s{\tilde \eta }^{\psi
}(\zeta ,z,0),\mbox{ }^s\zeta ) := [{\tilde {\partial }}_{\mbox{
}^s\zeta ,\mbox{ }^sz} \mbox{ }^s{\tilde {\eta }}^{\psi }(\zeta
,z,0)] \wedge
d\mbox{ }^s{\tilde {\zeta }}$, \\
${\check{\phi }}_{p,u}(\mbox{ }^s{\tilde \eta }^{\psi }(\zeta ,z,0),
\mbox{ }^s\zeta ):= (i_{2u} \{ [{\tilde {\partial }}_{\mbox{
}^s\zeta ,\mbox{ }^sz} \mbox{ }^s{\tilde {\eta }}^{\psi }(\zeta
,z,0)]i_{2u} \} ) \wedge (i_{2u}(d\mbox{ }^s\zeta i_{2u}))$ for each $u>0$,  \\
$(6)\quad \check{w}_{2^p}(\mbox{ }^s{\tilde \eta }^{\psi }(\zeta
,z,\lambda ),\mbox{ }^sz) := C_p \{ \check{\phi }_{p,0}(\mbox{
}^s{\tilde \eta }^{\psi }(\zeta ,z,\lambda ),\mbox{ }^sz)\wedge
\check{\phi }_{p,1}(\mbox{ }^s{\tilde \eta }^{\psi }(\zeta
,z,\lambda ),\mbox{ }^sz)\wedge ... $ \\
 $\wedge \check{\phi
}_{p,2^{p-1}-1}(\mbox{ }^s{\tilde \eta }^{\psi }(\zeta
,z,\lambda ),\mbox{ }^sz)\} _{q_0(2^{p-1})}$, \\
where $C_p=([(2^p-2)!]2(2^{p-1}-1))^{-1}$; \\
$(7)\quad \check{w}_{2^p,u}(\mbox{ }^s{\tilde \eta }^{\psi }(\zeta
,z,\lambda ),\mbox{ }^sz):= \{ \check{\phi }_{p,0}(\mbox{ }^s{\tilde
\eta }^{\psi }(\zeta ,z,\lambda ),\mbox{ }^sz)\wedge ... \wedge
\check{\phi }_{p,u-1}(\mbox{ }^s{\tilde \eta }^{\psi }(\zeta
,z,\lambda ),\mbox{ }^sz)\wedge $ \\
${\check{\phi }'}_{p,u} (\mbox{
}^s{\tilde \eta }^{\psi }(\zeta ,z,\lambda ),\mbox{ }^sz)\wedge
\check{\phi }_{p,u+1}(\mbox{ }^s{\tilde \eta }^{\psi }(\zeta
,z,\lambda ),\mbox{ }^sz)\wedge ... \wedge \check{\phi
}_{p,2^{p-1}-1}(\mbox{ }^s{\tilde \eta
}^{\psi }(\zeta,z,\lambda ),\mbox{ }^sz) \} _{q_0(2^{p-1})}$, \\
analogously to $(3-7)$ there are defined ${\check{\phi
}}_{p,u}(\mbox{ }^s {\tilde \psi }(\zeta ,z)),$ ${\check{\phi
}'}_{p,u}(\mbox{ }^s{\tilde \psi } (\zeta ,z))$,
$\check{w}_{2^p}(\mbox{ }^s{\tilde {\psi }}(\zeta ,z),\mbox{ }^sz)$,
$\check{w}_{2^p,u}(\mbox{ }^s{\tilde {\psi }}(\zeta ,z),\mbox{
}^sz)$ for each $u \ge 0$ with $\mbox{ }^s{\tilde \psi }(\zeta ,z)$
instead of $\mbox{ }^s{\tilde {\eta }}^{\psi } (\zeta ,z,\lambda );$
$$(8)\quad \phi _{\zeta ,z}:=\phi _{\zeta ,z}(\psi (\zeta ,z);\zeta ):=
{C'}_p <\psi (\zeta ,z); \zeta -z> ^{-2^{p-1}n} $$
$$\sum_{s=1}^n \sum_{u=0}^{2^{p-1}-1} \{ {\check{w}}_{2^p}( \mbox{ }^1{\tilde \psi }(\zeta
,z),\mbox{ }^1z) \wedge ... \wedge {\check{w}}_{2^p}(\mbox{
}^{s-1}{\tilde \psi }(\zeta ,z),\mbox{ }^{s-1}z)\wedge $$
$${\check{w}}_{2^p,u} (\mbox{ }^s{\tilde \psi }(\zeta ,z),\mbox{ }^sz)\wedge
{\check{w}}_{2^p}(\mbox{ }^{s+1}{\tilde \psi } (\zeta ,z),\mbox{
}^{s+1}z) \wedge ... \wedge {\check{w}}_{2^p}(\mbox{ }^n{\tilde \psi
}(\zeta ,z), \mbox{ }^n\zeta ) \} _{q_0(n)} ;$$
$$(9)\quad {\bar \phi }_{\zeta , z, \lambda }:={\bar \phi }_{\zeta , z, \lambda }
(\psi (\zeta ,z); \zeta ):= $$
$${C'}_p \sum_{s=1}^n \sum_{u=0}^{2^{p-1}-1} \{ {\check{w}}_{2^p}(
\mbox{ }^1{\tilde {\eta }}^{\psi }(\zeta ,z,\lambda ),\mbox{ }^1z)
\wedge ... \wedge {\check{w}}_{2^p}(\mbox{ }^{s-1}{\tilde {\eta
}}^{\psi }(\zeta ,z,\lambda ),\mbox{ }^{s-1}z)\wedge $$
$${\check{w}}_{2^p,u} (\mbox{ }^s{\tilde {\eta }}^{\psi }
(\zeta ,z,\lambda ),\mbox{ }^sz)\wedge {\check{w}}_{2^p}(\mbox{
}^{s+1}{\tilde {\eta }}^{\psi } (\zeta ,z,\lambda ),\mbox{ }^{s+1}z)
\wedge ... \wedge {\check{w}}_{2^p}(\mbox{ }^n{\tilde {\eta }}^{\psi
}(\zeta ,z,\lambda ), \mbox{ }^n\zeta ) \} _{q_0(n)} .$$ If $f$ is a
bounded differential form on $U$, then define the integral
operators:
$$(10)\quad (L^{\psi }_{\partial U}f)(z):=\int_{\zeta \in \partial U}
f(\zeta ) \wedge \phi _{\zeta ,z}(\psi (\zeta ,z);\zeta ),$$
$$(11)\quad (R^{\psi }_{\partial U}f)(z):=\int_{\zeta \in \partial U,
0\le \lambda \le 1} f(\zeta ) \wedge {\bar \phi }_{\zeta , z,
\lambda }( \psi (\zeta ,z);\zeta ).$$
\par {\bf 3.5. Theorem.} {\it Let $U$ be an open subset in ${\cal
A}_p^n$, $2\le p\in \bf N$, with a piecewise $C^1$-boundary and let
$\psi $ be an ${\cal A}_p$-boundary distinguishing map for $U$.
Suppose that $f$ is a continuous mapping $f: cl(U)\to {\cal A}_p$
such that ${\tilde
\partial }f$ is also continuous on $U$ in the sence of distributions
and has a continuous extension on $cl (U)$. Then
$$(1)\quad f=(L^{\psi }_{\partial U}f)-(R^{\psi }_{\partial U}
{\tilde \partial }f)-(B_U{\tilde \partial }f)\mbox{ on }U,$$ where
the ${\cal A}_p$ integral operators $B_U$, $L^{\psi }_{\partial U}$
and $R^{\psi }_{\partial U}$ are given by Equations $3.1(13)$,
$3.4(10,11)$.}
\par {\bf Proof.} There is the decomposition:
\par $(2)\quad {\bar \phi }_{\zeta , z, \lambda }=\sum_{q=0}^{2^{p-1}n-1}
\Upsilon ^{\psi }_q (\zeta ,z,\lambda ),$ \\
where $\Upsilon ^{\psi }_q (\zeta ,z,\lambda )$ is a differential
form with all terms of degree $q$ by $z$ and $\tilde z$ and their
multiples on ${\cal A}_p$ constants and of degree $(2^pn-q-1)$ by
$(\zeta , \lambda )$ (including $\tilde \zeta $ and multiples of
$\zeta $ and $\tilde \zeta $ on ${\cal A}_p$ constants). A
differential form $f$ has Decomposition $3.2(12)$. If $\psi (z)$ is
an ${\cal A}_p$ $z$-superdifferentiable nonzero function on an open
set $V$ in ${\cal A}_p^n$, then differentiating the equality $(\psi
(z))(\psi (z))^{-1}=e$ gives $(\psi (z))\{ d_z(\psi (z))^{-1}.h \} =
- (d_z\psi (z).h)(\psi (z))^{-1}$
for each $z\in V$ and each $h\in {\cal A}_p^n$. Then \\
$\int_{\zeta \in \partial U, 0\le \lambda \le 1} f_r(\zeta )\wedge
\Upsilon ^{\psi }_q(\zeta ,z,\lambda )=0$
for each $r\ne q+1$, \\
since $dim (\partial U)=2^pn-1$, $d\lambda \wedge d\lambda =0$ and
$d\lambda $ commutes with each $b\in {\cal A}_p$. Therefore,\\
$(3)\quad R^{\psi }_{\partial U}f_r= \int_{\zeta \in \partial U,
0\le \lambda \le 1} f_r(\zeta )\wedge \Upsilon ^{\psi }_{r-1}(\zeta
,z,\lambda ) \mbox{ for each } 1\le r \le 2^{p-1}n$ and $R^{\psi
}_{\partial U}f_r=0$
for $r=0$ or $r>2^{p-1}n$. In particular, if $f=f_1$, then \\
$(4)\quad R^{\psi }_{\partial U}f_1= \int_{\zeta \in \partial U,
0\le \lambda \le 1} f_1(\zeta )\wedge {\bar \phi }_{\zeta , \lambda
} (\psi (\zeta ,z); \zeta ),$ \\
where ${\bar \phi }_{\zeta , \lambda } (\psi (\zeta ,z); \zeta ),$
is obtained from ${\bar \phi }_{\zeta , z, \lambda } (\psi (\zeta
,z); \zeta )$ by substituting all ${\tilde \partial }_{\mbox{
}^s\zeta ,\mbox{ }^sz}$ in Formulas $3.4(3-7,9,11)$ on ${\tilde
\partial }_{\mbox{ }^s\zeta }$.
On the other hand, with the help of Formulas $3.2(1-9,15)$ each
${\cal A}_p$ external derivative ${\tilde \partial }_{\mbox{
}^s\zeta }$ can be replaced on $d_{\mbox{ }^s\zeta }$ in ${\bar \phi
}_{\zeta , \lambda }(\psi (\zeta ,z); \zeta )$ in Formula $(4)$.
For $\phi _{\zeta ,z}$ there is the decomposition: \\
$(5)\quad \phi _{\zeta ,z }=\sum_{q=0}^{2^{p-1}n-1} \Upsilon ^{\psi
}_q (\zeta ,z),$ where $\Upsilon ^{\psi }_q(\zeta , z)$ is a
differential form with all terms of degree $q$ by $z$ and $\tilde z$
and their multiples on ${\cal A}_p$ constants and of degree
$2^pn-q-1$ by $\zeta $ and $\tilde \zeta $
and their multiples on ${\cal A}_p$ constants. Therefore, \\
$(6)\quad L^{\psi }_{\partial U}f_r=\int_{\zeta \in \partial U}
f_r\wedge \Upsilon ^{\psi }_r(\zeta ,z)$
for each $0\le r\le 2^{p-1}n-1$ \\
and $L^{\psi }_{\partial U}f_r=0$ for $r\ge 2^{p-1}n$.
In particular, for $f=f_0$: \\
$(7)\quad L^{\psi }_{\partial U}f_0= \int_{\zeta \in \partial U}
f_0(\zeta )\phi _{\zeta } (\psi (\zeta ,z); \zeta ),$ \\
where $\phi _{\zeta } (\psi (\zeta ,z)),$ is obtained from $\phi
_{\zeta , z} (\psi (\zeta ,z); \zeta )$ by substituting all ${\tilde
\partial }_{\zeta ,z}$ in Formulas $3.4.(3-10)$ on ${\tilde
\partial }_{\zeta }$.
\par In view of Formula $3.2(1)$ it remains to prove, that
$R^{\psi }_{\partial U}{\tilde \partial }f= L^{\psi }_{\partial
U}f-B_{\partial U}f$ on $U$. For each $\zeta $ in a neighborhood of
$\partial U$ there is the identity:
\par $(8)\quad <\eta ^{\psi }(\zeta , z, \lambda ); \zeta -z>=1$
for each $0\le \lambda \le 1$, hence $d_{\zeta ,z,\lambda } <\eta
^{\psi }(\zeta ,z,\lambda );\zeta -z>=0$. By Formulas $3.2(15-18):$
\par $(9)\quad d_{\zeta ,\lambda }{\bar \phi }_{\zeta , z, \lambda }=0.$
From Identities $3.2(15-18)$ it follows, that
\par $(10)\quad \partial _{\zeta }f\wedge {\bar \phi }_{\zeta , \lambda }=0.$
Therefore, from $(4), (9), (10)$ it follows, that
\par $(11)\quad d_{\zeta ,\lambda }[f(\zeta ){\bar \phi }_{
\zeta , \lambda }]=[{\tilde \partial }_{\zeta }f(\zeta )]\wedge
{\bar \phi }_{\zeta , \lambda },$ since ${\tilde \partial }_{\zeta
}f(\zeta )= \sum_{s=1}^n(\partial  f(\zeta , {\tilde \zeta
})/\partial \mbox{ }^s{\tilde \zeta }).d\mbox{ }^s{\tilde \zeta }.$
Due to Formulas $3.2(15-17)$ and $3.4(1-9)$:
\par $(12)\quad {\bar \phi }_{\zeta , \lambda }|_{\lambda =0}=
\phi _{\zeta }$, ${\bar \phi }_{\zeta , \lambda }|_{\lambda =1}=
\theta _z(\zeta )$. \\
If $\Upsilon (\zeta ,z,\lambda )$ is a differential form over ${\cal
A}_p$, then \\
$\Upsilon (\zeta ,z,\lambda ) = \sum_{s=0}^{2^p-1}\Psi _s(\zeta
_0,...,\zeta_{2^p-1},z_0,...,z_{2^p-1},\lambda )i_s,$ \\
where $\zeta =\zeta _0i_0+...+\zeta _{2^p-1}i_{2^p-1}$, $\zeta , z
\in {\cal A}_p$, $\zeta _0,...,\zeta _{2^p-1}, z_0,...,z_{2^p-1},
\lambda \in \bf R$, $\{ i_0,...,i_{2^p-1} \} $ denotes the set of
standard generators of ${\cal A}_p$, $\Psi _s$ is with values in
$\bf R$ for each $s=0,...,2^p-1$. From the Stoke's formula for
vector-valued differential forms, in particular, for \\
$[f(\zeta ){\bar \phi }_{\zeta ,z,\lambda } (\psi (\zeta ,z);\zeta
)]$ on $\partial U\times [0,1]$ and Formulas $(4), (7), (11), (12)$
above it follows the statement of this theorem.
\par {\bf 3.6. Corollary.} {\it Let conditions of Theorem 3.5
be satisfied and let $f$ be an ${\cal A}_p$ holomorphic function on
$U$, then $f=L^{\psi }_{\partial U}f$ on $U$.}
\par {\bf 3.7. Remark.} For $n=1$ Formula $3.2.(1)$
produces another analog of the Cauchy-Green formula (see Theorem
$2.1$ and Remark $2.1.1$) without using the ${\cal A}_p$ line
integrals. This is caused by the fact that the dimension of ${\cal
A}_p$ over $\bf R$ is greater, than $2$: $\quad dim_{\bf R}{\cal
A}_p=2^p$, that produces new integral relations. Theorem $3.2$ can
be used instead of Theorem $2.1$ to prove theorems $2.3$ and $2.11$
(with differential forms of Theorem $3.2$ instead of differential
forms of Theorem $2.1$). \par  If $\psi (\zeta ,z)=\zeta -z$, then
$L^{\psi }_{\partial U}= B_{\partial U}$ and $R^{\psi }_{\partial
U}=0$, hence Formula $3.5.(1)$ reduces to Formula $3.2.(1)$. For a
function $f$ or a $1$-form $f$ Formulas $3.2(13,14)$ respectively
are valid as well for ${\bar \theta }(\zeta ,z)$ instead of $\theta
_z(\zeta )$, where $d_{\zeta ,z}{\bar \theta } (\zeta ,z)=0$ for
each $\zeta \ne z$. A choice of $\phi _{p,s}$ and $w_{2^p}$ is not
unique, for example, $d{\tilde \zeta }\wedge d\zeta \wedge d\zeta
\wedge d\zeta $ may be taken in $\bf H$, since it gives up to a
multiplier $Ce$, where $C$ is a real constant, the canonical volume
element in $\bf H$ and $d\zeta \wedge d\zeta \wedge d\zeta \wedge
d\zeta =0.$ \par  Formulas $3.2.(1)$ and $3.5.(1)$ for functions of
${\cal A}_p$ variables are the ${\cal A}_p$ analogs of the
Martinelli-Bochner and the Leray formulas for functions of complex
variables respectively, where $\psi (\zeta ,z)$ is the ${\cal A}_p$
analog of the Leray complex map (see \S 3.4). In the ${\cal A}_p$
case the algebra of differential forms bears the additional
gradation structure and have another properties, than in the complex
case (see also \S \S 2.8 and 3.7 \cite{luoyst}). Lemma $3.9$ below
shows, that the ${\cal A}_p$ boundary distinguishing maps exist.
\par {\bf 3.8. Definitions and Notations.}
Let a subset $U$ in ${\cal A}_p^n$, $2\le p\in \bf N$ or $p=\Lambda
$, be given by the equation:
\par $(1)$ $U:=\{ z\in {\cal A}_p^n:$ $\rho (z)<0 \} $,
where $\rho $ is a real-valued $C^2$-function such that there exists
a constant $\epsilon _0>0$ for which:
\par $(2)$ $\sum_{l,m=1}^{2^pn}(\partial ^2\rho (z)/
\partial x_l\partial x_m)t_lt_m\ge \epsilon _0|t|^2$
for each $t\in \bf R^{2^pn}$ for finite $p$ with \\
$\mbox{ }^lz= \sum_{m=1}^{2^p}x_{2^p(l-1)+m}S_m$, $S_m := i_{m-1}$
for each $m$, $x_l\in \bf R$; or
\par $(2)'$ $\sum_{l,m=1}^n(\partial ^2\rho (z)/
\partial \mbox{ }^lz\partial \mbox { }^m{\tilde z}).(\mbox{ }^lh,\mbox{ }^m{\tilde h})\ge
\epsilon _0 \| h \|^2$ for each $h\in {\cal A}_p^n$ for infinite
$p=\Lambda $, where $z=(\mbox{ }^1z,..., \mbox{ }^nz)$, $\mbox{
}^lz\in {\cal A}_p$. Then $U$ is called strictly convex open subset
(with $C^2$-boundary). Let
\par $(3)$ $w_{\rho } (z):=(\partial \rho (z)/ \partial \mbox{ }^1z,...,
\partial \rho (z)/ \partial \mbox{ }^nz)$.
\par Put for finite $2\le p\in \bf N$
\par $v_{\rho }(z):=\sum_{m=1}^{2^p}(w_{\rho }.S_m)S_m,$ where as usually
$w_{\rho }.S_m=(d_z\rho (z)).S_m$ is the differential of the
function $\rho $.
\par For infinite $p=\Lambda $ take a collar neighborhood $V$ for
$\partial U$ such that for each $\zeta \in V$ there exists a unique
point $\xi \in \partial U$, for which $\zeta $ belongs to the
segment of the straight line intersecting $\partial U$ at the point
$\xi $ along an outside normal (perpendicular) vector $n_{\xi }$ to
$\partial U$ at the point $\xi \in
\partial U$, $\xi =\xi (\zeta )$. Put \par $<v_{\rho }(\zeta
);h>:=\sum_m <((\partial \rho (\zeta )/\partial \zeta ).S_m)S_m;h>$, \\
where $h\in {\cal A}_p^n$, that defines $v_{\rho }$, since $\| h \|
<\infty $ for each $h\in {\cal A}_p^n$ and $(\partial \rho (\zeta
)/\partial \zeta )\in C^1$ is a bounded operator for each $\zeta $.
\par {\bf 3.9. Lemma.} {\it Let the function $v_{\rho }$ be as in \S 3.8.
Then $v_{\rho }$ is the ${\cal A}_p$-boundary distinguishing map for
$U$.}
\par {\bf Proof.} Since $S_mS_l=(-1)^{\kappa (S_m) \kappa (S_l)}S_lS_m$
for each $m\ne l$,
where $\kappa (S_1)=0$, $\kappa (S_m)=1$ for each $m>1$, then \\
$<v_{\rho }(\zeta );\zeta -z>+<\zeta -z;v_{\rho }(\zeta )>=2
\sum_{l=1}^{2^pn}(\partial \rho (\zeta )/\partial x_l) x_l(\zeta -z)$, \\
when $p$ is finite, where $x_l=x_l(\zeta )$ and $x_l(\zeta -z)$ are
real coordinates corresponding to $\zeta $ and $\zeta -z$. For
infinite $p=\Lambda $ there is the equality: \\
$<v_{\rho }(\zeta );\zeta -z> + <\zeta -z;v_{\rho }(\zeta )>=2 Re [
\sum_{l=1}^n (\partial \rho (\zeta )/\partial \mbox{ }^l\zeta ).
(\mbox{ }^l\zeta -\mbox{ }^lz)]$. \\
By the Taylor's theorem: $\rho (z)= \rho (\zeta ) - <v_{\rho }(\zeta
);\zeta -z>/2 - <\zeta -z; v_{\rho }(\zeta )>/2 + \sum_{l,m=1}^n
(\partial ^2\rho (\zeta )/ \partial \mbox{ }^l\zeta \partial \mbox{
}^m\zeta ).((\mbox{ }^l\zeta - \mbox{ }^l z),(\mbox{ }^m\zeta
-\mbox{ }^m z)/2 +o(|\zeta -z|^2)$. Therefore, there exists a
neighborhood $V$ of $\partial U$ and $\epsilon _1>0$ such that \\
$(1)\quad (<v_{\rho }(\zeta );\zeta -z>+<\zeta -z; v_{\rho }(\zeta
)>)/2 \ge \rho (\zeta )-\rho (z)+\epsilon _0 |\zeta -z|^2/4$ for
each $\zeta \in V$ and $|\zeta -z|\le \epsilon _1$, where
$a=\sum_ma_mS_m$ for each $a\in {\cal A}_p$, $a_m$ are reals. If
$z\in U$, $\zeta \in
\partial U$, $|z-\zeta |\le \epsilon _1$, then by $(1)$: $(<v_{\rho
}(\zeta ); \zeta -z>+<\zeta -z; v_{\rho }(\zeta )>) \ge -\rho
(z)>0$. If $|\zeta -z|>\epsilon _1$, put $z_1:=(1-\epsilon _1|\zeta
-z|^{-1})\zeta + \epsilon _1|\zeta -z|^{-1}z$, then $\zeta
-z_1=\epsilon _1|\zeta -z|^{-1} (\zeta -z)$, consequently,
$(<v_{\rho }(\zeta ); \zeta -z>+<\zeta -z; v_{\rho }(\zeta )>)_e/2
\ge -\rho (z_1)$. Evidently, $U$ is convex and $z_1\in U$.
\par {\bf 3.10. Theorem.} {\it Let $U$ be a strictly
convex open subset in ${\cal A}_p^n$, $2\le p\in \bf N$ or
$p=\Lambda $, (see $3.8.(1)$) and let $f$ be a continuous function
on $U$ with continuous ${\tilde \partial }f$ on $U$ in the sense of
distributions having a continuous extension on $cl (U)$ such that
$2.11.(i)$ is satisfied. Then there exists a function $u$ on $U$
which is a solution of the ${\tilde \partial }$-equation
$2.11.(ii)$.}
\par {\bf Proof.} Reduce the proof to the finite case as in \S 2.3.
In proofs of Theorems $2.3$ and $2.11$ take in Formula $3.5(1)$
$\chi f$ instead of $f$, which is possible due to Lemma $3.9$,
choosing $\psi = v_{\rho }$ and $supp (\chi )$ as a proper subset of
$U$. Then $L^{\psi }_{\partial U}\chi f=0$ and $R^{\psi }_{\partial
U}\chi f=0$, hence $\chi f=-B_U{\tilde
\partial }\chi f$. For each fixed $z\in U$ a subset $\mbox{ }^lU_{\eta }:=
\{ \xi \in {\cal A}_p: \rho (\mbox{ }^1z,..., \mbox{ }^{l-1}z,\xi
,\mbox{ }^{l+1}z,...,\mbox{ }^nz) <0 \} $ is strictly convex in
${\cal A}_p$ due to $3.8(1,2)$, where $\eta :=(\mbox{
}^1z,...,\mbox{ }^{l-1}z,\mbox{ }^{l+1}z,...,\mbox{ }^nz).$ Apply
$3.5(1)$ by a variable $\xi $ in $\mbox{ }^lU_{\eta }$, in
particular, for $l=1$, for which $v_{\rho }$ by the variable $\xi $
is the ${\cal A}_p$-boundary distinguishing map for $\mbox{
}^1U_{\eta }$. Therefore, $u(z):=-B_{\mbox{ }^1U_{\eta }}\mbox{
}^1{\hat f}(\xi ,\eta ). d{\tilde \xi }$ with $z=(\xi ,\eta )$, $\xi
\in \mbox{ }^1U_{\eta }$ solves the problem.
\section{Manifolds over Cayley-Dickson algebras}
\par {\bf 4.1. Definitions and Notations.}
Suppose that $M$ is an ${\cal A}_p$ manifold and let $RL(N,{\cal
A}_p)$ be the family of all right ${\cal A}_p$-superlinear operators
$A: {\cal A}_p^N\to {\cal A}_p^N$, where $2\le p\in \bf N$ or
$p=\Lambda $. Then an ${\cal A}_p$ holomorphic vector bundle $Q$ of
${\cal A}_p$ dimension $N$ over $M$ is a $C^{\infty }$-vector bundle
$Q$ over $M$ with the characteristic fibre ${\cal A}_p^N$ together
with an ${\cal A}_p$ holomorphic atlas of local trivializations:
$g_{a,b}: U_a\cap U_b\to RL(N,{\cal A}_p)$, where $U_a\cap U_b\ne
\emptyset $, $ \{ (U_a,h_a): a \in \Upsilon \} =At (Q)$, $\bigcup_a
U_a=M$, $U_a$ is open in $M$, $h_a: Q|_{U_a}\to U_a\times {\cal
A}_p^N$ is the bundle isomorphism, $(z,g_{a,b}(z)v)=h_a\circ
h_b^{-1}(z,v)$, $z\in U_a\cap U_b$, $v\in {\cal A}_p^N$. If $Y_N$ is
the underlying to ${\cal A}_p^N$ real vector space, then suppose,
that each $g_{a,b}$ induces a $\bf R$-linear isomorphism of $Y_N$
onto itself. Since $M$ has the real underlying manifold $M_{\bf R}$,
then there exists the tangent bundle $TM$ such that $T_xM$ is
isomorphic with ${\cal A}_p^n$ for each $x\in M$, since
$TU_a=U_a\times {\cal A}_p^n$ for each $a$, where $dim_{{\cal
A}_p}M=n$ is the ${\cal A}_p$ dimension of $M$. \par  If $X$ is a
Banach space over ${\cal A}_p$ (with left and right distributivity
laws relative to multiplications of vectors in $X$ on scalars from
${\cal A}_p$), then denote by $X^*_q$ the space of all additive $\bf
R$-homogeneous functionals on $X$ with values in ${\cal A}_p$.
Clearly $X^*_q$ is the Banach space over ${\cal A}_p$. Then $T^*M$
with fibres $(({\cal A}_p^n)_q)^*$ denotes the ${\cal A}_p$
cotangent bundle of $M$ and $\Lambda ^rT^*M$ denotes the vector
bundle whose sections are ${\cal A}_p$ $r$-forms on $M$, where
$S_bdx_b\wedge S_adx_a=-(-1)^{\kappa (S_a) \kappa (S_b)}
S_adx_a\wedge S_bdx_b$ for each $S_a\ne S_b\in \{ i_0,...,i_{2^p-1}
\} $, $dz=\sum_{m=1}^{2^p}dx_m S_m$, $z\in {\cal A}_p$, $x_b\in \bf
R$.
\par The ${\cal A}_p$ holomorphic Cousin data in $Q$ is a family
$ \{ f_{a,b}: a, b \in \Upsilon \} $ of ${\cal A}_p$ holomorphic
sections $f_{a,b}: U_a\cap U_b\to Q$ such that
$f_{a,b}+f_{b,l}=f_{a,l}$ in $U_a\cap U_b\cap U_l$ for each $a, b, l
\in \Upsilon $. A finding of a family $\{ f_a: a\in \Upsilon \} $ of
${\cal A}_p$ holomorphic sections $f_a: U_a\to Q$ such that
$f_{a,b}=f_a-f_b$ in $U_a\cap U_b$ for each $a, b \in \Upsilon $
will be called the ${\cal A}_p$ Cousin problem.
\par {\bf 4.2. Theorem.} {\it Let $M$ be an ${\cal A}_p$ manifold and $Q$ be an
${\cal A}_p$ holomorphic vector bundle on $M$, where $2\le p\in \bf
N$ or $p=\Lambda $. Then Conditions $(i,ii)$ are equivalent:
\par $(i)$ each ${\cal A}_p$ holomorphic Cousin problem in $M$
has a solution;
\par $(ii)$ for each ${\cal A}_p$ holomorphic section $f$ of $Q$
such that ${\tilde \partial }f=0$ on $M$, there exists a $C^{\infty
}$-section $U$ of $Q$ such that $(\partial u/\partial {\tilde
z})={\hat f}$ on $M$.}
\par {\bf Proof.} $(i)\rightarrow (ii)$. In view of Theorems $2.11$
and $3.10$ there exists an (open) covering $ \{ U_a : a \} $ of $M$
and $C^{\infty }$-sections $u_b: U_b\to Q$ such that $(\partial
u_b/\partial {\tilde z})={\hat f}$ in $U_b$. Then $(u_b-u_l)$ is
${\cal A}_p$ holomorphic in $U_b\cap U_l$ and their family forms the
${\cal A}_p$ holomorphic Cousin data in $Q$. Put $u:=u_b-h_b$ on
$U_b$, where $u_b-u_l=h_b-h_l$, $h_b: U_b\to Q$ is an ${\cal A}_p$
holomorphic section given by $(i)$.
\par $(ii)\rightarrow (i)$.
Take a $C^{\infty }$-partition of unity $\{ \chi _b: b \} $
subordinated to $ \{ U_b: b \} $ and $c_b:=-\sum_a\chi _af_{a,b}$ on
$U_b$, then $f_{l,b}=\sum_a\chi _a(f_{l,a}+f_{a,b})=c_l-c_b$ in
$U_l\cap U_b$, hence $(\partial c_l/\partial {\tilde z}) =(\partial
c_b/\partial {\tilde z})$ in $U_l\cap U_b$. By $(ii)$ there exists a
$C^{\infty }$-section $u: M\to Q$ with $(\partial u/\partial {\tilde
z})=(\partial c_b/\partial {\tilde z})$ on $U_b$ and $h_b:=c_b-u$ on
$U_b$ gives the solution.
\par {\bf 4.3. Definitions.}
Suppose $U$ is an open subset in ${\cal A}_p$, then a $C^2$-function
$\rho : U\to \bf R$ is called subharmonic (strictly subharmonic) in
$U$ if $\sum_{m=1}^{2^p}\partial ^2\rho /\partial x_m^2\ge 0$
($\sum_{m=1}^{2^p}\partial ^2\rho /\partial x_m^2>0$
correspondingly) for finite $p\ge 2$; or
\par $ (\partial ^2 \rho (z)
/\partial z\partial {\tilde z}).(\xi ,{\tilde \xi })\ge 0$ (or $>0$)
for each $z\in U$ and each $0\ne \xi \in {\cal A}_p$ for $p=\Lambda
$, where $z=\sum_{m=1}^{2^p}x_mS_m\in U$, where $x_m\in \bf R$ for
each $m$.
\par If $U$ is an open subset in
${\cal A}_p^n$, then a $C^2$-function $\rho : U\to \bf R$ such that
the function $\zeta \mapsto \rho (v+\zeta w)$ is subharmonic
(strictly subharmonic) on its domain for each $v, w\in {\cal A}_p^n$
is called plurisubharmonic (strictly plurisubharmonic
correspondingly) function, where $\zeta \in {\cal A}_p$.
\par A $C^v$-function $\rho $ on an ${\cal A}_p$ manifold $M$
is called a strictly plurisubharmonic exhausting $C^v$-function for
$M$, $2\le v\in \bf N$, if $\rho $ is a strictly plurisubharmonic
$C^v$-function on $M$ and for each $\alpha \in \bf R$ the set $ \{
z\in M: \rho (z)<\alpha \} $ is bounded in $M$.
\par {\bf 4.4. Theorem.} {\it Let $M$ be an ${\cal A}_p$ manifold
with strictly plurisubharmonic exhausting function $\rho $ such that
$\rho $ is a $C^{\omega }_{z,\tilde z}$-function and let $Q$ be an
${\cal A}_p$ holomorphic vector bundle on $M$, $U_{\alpha }:= \{
z\in M: \rho (z)<\alpha \} $ for $\alpha \in \bf R$, where $2\le
p\in \bf N$ or $p=\Lambda $.
\par $(i).$ Suppose that $d\rho (z)\ne 0$ for each $z\in \partial
U_{\alpha }$ for a marked $\alpha \in \bf R$. Then every continuous
section $f: cl (U_{\alpha })\to Q$ which is ${\cal A}_p$ holomorphic
on $U_{\alpha }$ can be approximated (uniformly for finite $p$) on
$cl (U_{\alpha })$ by ${\cal A}_p$ holomorphic sections of $Q$ on
$M$.
\par $(ii).$ For each continuous mapping $f: M\to Q$ such that
${\tilde \partial f}=0$ on $M$ there exists a continuous mapping $u:
M\to Q$ such that $\partial u/\partial {\tilde z}={\hat f}$ on $M$.}
\par {\bf Proof.} For a $C^{\omega }_{z,\tilde z}$-function
$\rho : U\to \bf R$ (that is, $\rho $ is locally analytic in
variables $(z,{\tilde z})$, ${\bf R}={\bf R}e\hookrightarrow {\cal
A}_p$) there is the identity: \\
$\sum_{l,m,a,b} (\partial ^2\rho /\partial \mbox{ }^lx_a\partial
\mbox{ }^mx_b)t_{2^p(l-1)+a}t_{2^p(l-1)+b}=$ \\
$\sum_{m,l}(\partial ^2 \rho (z)/\partial \mbox{ }^lz\partial \mbox{
}^mz).((\partial \mbox{ }^lz/\partial \mbox{ }^lx_a)t_{2^p(l-1)+a},
(\partial \mbox{ }^mz/\partial \mbox{ }^lx_b)t_{2^p(l-1)+b})$ \\
$= \sum_{m,l=1}^n (\partial ^2 \rho (z)/\partial \mbox{ }^lz\partial
\mbox{ }^m{\tilde z}).(\mbox{ }^l\xi ,\mbox{ }^m{\tilde \xi })$ for
finite $p$, since $\partial \mbox{ }^lz/\partial \mbox{ }^lx_a=S_a$,
$\partial \mbox{ }^l{\tilde z}/\partial \mbox{ }^lx_a= (-1)^{\kappa
(S_a)}S_a$, where $\mbox{ }^l\xi =\sum_{m=1}^{2^p}
t_{2^p(l-1)+m}S_m$, $S_m=i_{m-1}$ for each $m$, $\mbox{ }^lz=
\sum_{m=1}^{2^p} \mbox{ }^lx_mS_m$, $t_b\in \bf R$, $\mbox{ }^lx_m
\in \bf R$.  Therefore, a $C^{\omega }_{z,\tilde z}$-function $\rho
$ is strictly plurisubharmonic on $U$ if and only if
\par $(1)\quad \sum_{m,l=1}^n (\partial ^2 \rho (z)
/\partial \mbox{ }^lz\partial \mbox{ }^m{\tilde z}).(\mbox{ }^l\xi
,\mbox{ }^m{\tilde \xi })>0$ for each $z\in U$ and each $0\ne \xi
\in {\cal A}_p^n$, where $\xi =(\mbox{ }^1\xi ,...,\mbox{ }^n\xi )$
(see also \S 2 \cite{luoyst,luoystoc}). Consider a proper bounded
closed subset $A$ in $M$ such that $d\rho (z)\ne 0$ for each $z\in
A$. Then for each $\epsilon >0$ there exists a strictly
plurisubharmonic function $\rho _{\epsilon }: M\to \bf R$ such that
$\rho e$ is a $C^{\omega }_{z,\tilde z}$-function on $M$ and
$(i-iii)$ are fulfilled:
\par $(i)$ $\rho -\rho _{\epsilon }$ together with its first
and second derivatives is not greater than $\epsilon $ on $M$;
\par $(ii)$ the set $Crit (\rho _{\epsilon }):=\{ z\in M:
d\rho _{\epsilon }(z)=0 \} $ is discrete in $M$;
\par $(iii)$ $\rho _{\epsilon }=\rho $ on $A$ (see also Lemma
$2.1.2.2$ \cite{henlei} in the complex case).
\par The space $C^{\omega }_z(U,{\cal A}_p)$ is dense in
$C^0(U,{\cal A}_p)$ for each open $U$ in ${\cal A}_p^n$ (see \S 2.7
and Theorem $3.28$ in \cite{luoyst,luoystoc}). Suppose $\beta \in
\bf R$ and $d\rho (z)\ne 0$ for $z\in \partial U_{\beta }$ and $f:
cl(U_{\beta })\to Q$ is a continuous section ${\cal A}_p$
holomorphic on $U_{\beta }$. Therefore, for each $\beta \le \alpha
<\infty $ if $d\rho (z)\ne 0$ for each $z\in \partial U_{\alpha }$,
then $f$ can be approximated (uniformly for finite $p$) on
$cl(U_{\beta })$ by continuous sections on $cl(U_{\alpha })$ that
are holomorphic on $U_{\alpha }$. There exists a sequence $\beta
<\alpha _1<\alpha _2<...$ such that $\lim_l \alpha _l=\infty $ and
$d\rho (z)\ne 0$ for each $z\in
\partial U_{\alpha _l}$, since $Crit (\rho )$ is discrete. For each
$\epsilon >0$ and each natural number $s$ satisfying $2\le s\le p$
there exists a continuous section $f_l: cl(U_{\alpha _l})\to Q$ such
that $f_l$ is ${\cal A}_p$ holomorphic on $U_{\alpha _l}$ and $\|
f_{l+1}-f_l \|_{C^0(U_{\alpha _l,s})}< \epsilon 2^{-l-1}$ for each
$l\in \bf N$, where $f_0:=f$, $M_s$ denotes the ($2^s$-dimensional
over $\bf R$) closed submanifold in $M$ induced by the embedding of
${\cal A}_s$ into ${\cal A}_p$, when $p$ is infinite, or put $M_s=M$
for finite $p$ taking $s=p$, $U_{\alpha _l,s}:=U_{\alpha _l}\cap
M_s$. Therefore, the sequence $ \{ f_l: l\in {\bf N} \} $ converges
to the ${\cal A}_s$ holomorphic section $g: M_s\to Q$ uniformly on
each compact subset $P$ in $M_s$ and $\| f-g \|_{C^0(U_{\beta ,s})}<
\epsilon $.
\par The second statement $(ii)$ follows from $(i)$ and
Theorems $2.11, 3.10$, since $Crit (\rho )$ is discrete in $M$ and
there exists a sequence of continuous $Q$-valued functions on $cl
(U_{\alpha _l})$ such that $\partial u_l/\partial {\tilde z}= \hat
f$ in the sense of distributions on $U_{\alpha _l}$,
$\bigcup_lU_{\alpha _l}=M$ (see also the complex case in \S 2.12.3
\cite{henlei} mentioning, that Lemma $2.12.4$ there can be
reformulated and proved for an ${\cal A}_p$ manifold $M$ on ${\cal
A}_p^n$ instead of a complex manifold on $\bf C^n$).
\par {\bf 4.5. Definitions.} Let $M$ be an ${\cal A}_p$ manifold
(see \S 2.10), where $2\le p\in \bf N$ or $p=\Lambda $. For a closed
bounded subset $G$ in $M$ put: ${\hat G}^{\cal H}_M:=\{ z\in M:
|f(z)|\le \sup_{\zeta \in G} |{\hat f}(\zeta )| \quad \forall f\in
{\cal H}(M) \} $. Such ${\hat G}^{\cal H}_M$ is called the ${\cal
H}(M)$-hull of $G$. If $G={\hat G}^{\cal H}_M$, then $G$ is called
${\cal H}(M)$-convex. An ${\cal A}_p$ manifold $M$ is called ${\cal
A}_p$ holomorphically convex if for each closed bounded subset $G$
in $M$ the set ${\hat G}^{\cal H}_M$ is closed and bounded.
\par An ${\cal A}_p$ manifold $M$ with a countable atlas $At(M)$
having dimension $n$ over ${\cal A}_p$ and satisfying $(i,ii)$:
\par $(i)$ $M$ is ${\cal A}_p$ holomorphically convex;
\par $(ii)$ for each $z\in M$ there are $\mbox{ }^1f,...,\mbox{ }^nf
\in {\cal H}(M)$ and there exists a neighbourhood $U$ of $z$ such
that the map $U\ni \zeta \mapsto (\mbox{ }^1f(\zeta ),...,\mbox{
}^nf (\zeta ))$ is ${\cal A}_p$ biholomorphic (see \S 2.6), then $M$
is called an ${\cal A}_p$ Stein manifold.
\par {\bf 4.6. Remark.} If $M_1$ and $M_2$ are two ${\cal A}_p$
Stein manifolds, then $M_1\times M_2$ is an ${\cal A}_p$ Stein
manifold. If $N$ is a closed ${\cal A}_p$ submanifold of a ${\cal
A}_p$ Stein manifold $M$, then $N$ is also an ${\cal A}_p$ Stein
manifold.
\par {\bf 4.7. Theorem.} {\it Let $M$ be an ${\cal A}_p$ Stein manifold,
where $2\le p\in \bf N$ or $p=\Lambda $. Then for each ${\cal
H}(M)$-convex closed bounded subset $P$ in $M$, $P\ne M$ and each
neighborhood $V_P$ of $P$ there exists a strictly plurisubharmonic
exhausting $C^{\omega }_{z,{\tilde z}}$-function $\rho $ on $M$ such
that $\rho <0$ on $P$ and $\rho
>0$ on $M\setminus V_P$.}
\par The {\bf proof} of this theorem is analogous to that of
Theorem $2.3.14$ \cite{henlei} in the complex case taking $\rho
(z):=-1+\sum_{l=1}^{\infty }\sum_{k=1}^{N(l)}
f^k_l(z)(f^k_l(z))^{\tilde .}$ for each $z\in M$, where $f^k_l\in
C^{\omega }_{z,\tilde z}$, $M=\bigcup_lP_l$, $P_l\subset Int
(P_{l+1})$ for each $l\in \bf N$, each $P_l$ is ${\cal
H}(M)$-convex, $\sum_{k=1}^{N(l)}|f^k_l(z)|^2<2^{-l}$ for each $z\in
P_l$, $\sum_{k=1}^{N(l)}|f^k_l(z)|^2>l$ for each $z\in
P_{l+2}\setminus U_l$, $U_l:=Int (P_{l+1})$, with the rank $rank
[(\partial f^k_l/\partial \mbox{
}^mz)^{k=1,...,N(l)}_{m=1,...,n}]=2^pn$ over $\bf R$ for each $z\in
P_l$ for finite $p$ or $(f^1,...,f^{N(l)})(z)$ is regular for
infinite $p$ (see Definitions 2.6 and 4.3).
\par {\bf 4.8. Theorem.} {\it An ${\cal A}_p$ manifold $M$ is an
${\cal A}_p$ Stein  manifold if and only if there exists a strictly
plurisubharmonic exhausting $C^{\omega }_{z,\tilde z}$-function
$\rho $ on $M$, then $ \{ z\in M: \rho (z)\le \alpha \} $ is ${\cal
H}(M)$-convex for each $\alpha \in \bf R$, where $2\le p\in \bf N$
or $p=\Lambda $.}
\par {\bf Proof.} The necessity follows from Theorem $4.7$.
To prove sufficiency suppose $\eta =(\mbox{ }^1\eta ,..., \mbox{
}^n\eta )$ are ${\cal A}_p$ holomorphic coordinates in a
neighbourhood $V_{\xi }$ of $\xi \in M$. Consider
\par $(1)\quad u(z):=2\sum_{l,m=1}^n<v_{\rho }(\xi );\eta (z)-\eta (\xi )>$ \\
$+\sum_{l,m=1}^n(\partial ^2\rho (\xi )/\partial \mbox{ }^l\eta
\partial \mbox{ }^m\eta ).[(\mbox{ }^l\eta (z)-\mbox{ }^l\eta (\xi )),
(\mbox{ }^m\eta (z)-\mbox{ }^m\eta (\xi ))]e$, \\
where $v_{\rho }(\xi )$ is given by $3.8.(3)$. Then $u$ is
holomorphic in $V_{\xi }$ and $u(\xi )=0$. By Lemma $3.9$:
\par $(2)\quad (u(z)+{\tilde u}(z))/2=\rho (z)-\rho (\xi )-
\sum_{l,m=1}^n (\partial ^2\rho (\xi )/\partial \mbox{ }^l\eta
\partial \mbox{ }^m{\tilde {\eta }}).[(\mbox{ }^l\eta (z)- \mbox{ }^l\eta (\xi )),
(\mbox{ }^m\eta (z)- \mbox{ }^m\eta (\xi ))^{\tilde .}] +o(|\eta
(\xi )-\eta (z)|^2)$. From the strict plurisubharmonicity of $\rho $
it follows, that there exists $\beta >0$ and $V_{\xi }$ such that
\par $(3)\quad (u(z)+{\tilde u}(z))/2< \rho (z)-\rho (\xi )
-\beta |\eta (z)-\eta (\xi )|^2$ for each $z\in V_{\xi }$. Then
$\exp (u(\xi ))=1$ and $|\exp (u(z))|<1$ for each $\xi \ne z\in
cl(U_{\alpha })\cap V_{\xi }$ (see Corollary $3.3$ \cite{luoyst},
Corollary 3.3 and Note 3.6.3 \cite{luoystoc}).
\par If $g: {\bf R}\to {\cal A}_p$ is a $C^{\infty }$-function
with bounded (closed) support, then $g(z{\tilde z})=:\chi (z)$ is a
$C^{\infty }$-function on ${\cal A}_p^n$ with bounded closed support
such that $\chi $ is ${\cal A}_p$ $(z,{\tilde
z})$-superdifferentiable. Therefore, there exists a neighbourhood
$W_{\xi }\subset V_{\xi }$ of $\xi $ and an infinitely $(z,{\tilde
z})$-superdifferentiable function $\chi $ such that $\chi |_{W_{\xi
}}=1$, $supp (\chi )$ is a proper subset of $V_{\xi }$,
consequently, \par $\lim_{m\to \infty } \| \exp (mu(z)) (\partial
\chi (z)/\partial
{\tilde z}) \| _{C^0(U_{\alpha })}=0$, \\
where $(\partial \chi
(z)/\partial {\tilde z})=(\partial \chi (z)/
\partial \mbox{ }^1{\tilde z},...,\partial \chi (z)/\partial \mbox{ }^n
{\tilde z}).$ In view of Theorem $3.10$ there exist continuous
functions $v_m$ on $cl (U_{\alpha })$ such that \\
$(\partial
v_m/\partial {\tilde z})= \exp (mu(z))(\partial \chi /\partial
{\tilde z})$ in $U_{\alpha }$ and $\lim_{m\to \infty } \| v_m
\|_{C^0(U_{\alpha })}=0$.
\par Put $g_m(z):=\exp (mu(z))\chi (z)-v_m(z)+v_m(\xi )$, hence
$g_m$ is continuous on $cl(U_{\alpha })$ and ${\cal A}_p$
holomorphic on $U_{\alpha }$. Since $supp (\chi )$ is the proper
subset in $V_{\xi }$, then $g_m(\xi )=1$ for each $m\in \bf N$,
$sup_m \| g_m \|_{C^0(U_{\alpha })}<\infty $ and for each bounded
closed subset $P$ in $cl(U_{\alpha })\setminus \{ \xi \} $ there
exists $\lim_m \| g_m \|_{C^0(P)}=0$. In view of Theorem $4.4.(ii)$
there exists a sequence of functions $f_m\in {\cal H}(M)$ and
$C=const<\infty $ such that $(a)$ $f_m(\xi )=1$ for each $m\in \bf
N$; $(b)$ $\| f_m \|_{C^0(U_{\alpha })}\le C$ for each $m\in \bf N$;
$(c)$ $\lim_{m\to \infty } \| g_m \|_{C^0(P)}=0$ for each closed
bounded subset $P\subset cl(U_{\alpha })\setminus \{ \xi \} $.
\par Consider an ${\cal A}_p$ holomorphic function $f$ on a neighborhood
of $\xi $ such that $f(\xi )=0$. Put $\phi _m:=f\exp (mu)\partial
\chi /\partial {\tilde z}$, then $supp (\phi _m)$ is the proper
subset in $V_{\xi }\setminus W_{\xi }$. In view of Inequality $(3)$
there exists $\delta >0$ such that $\lim_m \| \phi _m
\|_{C^0(U_{\alpha +\delta })}=0$. As in \S 4.4 it is possible to
assume, that $Crit (\rho )$ is discrete in $M$. Take $0<\epsilon
<\delta $ such that $d\rho \ne 0$ on $\partial U_{\alpha +\epsilon
}$. In view of Theorem $4.4.(ii)$ there exists a continuous function
$v_m$ on $cl(U_{\alpha +\epsilon })$ such that $\partial
v_m/\partial {\tilde z}={\hat \phi }_m$ on $U_{\alpha +\epsilon }$
and $\lim_m \| v_m \|_{C^0(U_{\alpha +\epsilon })} =0$. Each $v_m$
is ${\cal A}_p$-holomorphic on $W_{\xi }$, since $\phi _m=0$ on
$W_{\xi }$, hence $\lim_m\partial v_m(\xi )=0$. Since $f(\xi )=u(\xi
)=0$ and $\chi =1$ on $W_{\xi }$, then $\partial g_k(\xi )/\partial
\xi =\partial f(\xi )/\partial \xi -\partial v_k(\xi )/\partial \xi
$, where $g_k:=f(\chi \exp (ku))-v_k$. In view of Theorem $4.4.(i)$
there exists $f_m\in {\cal H}(M)$ such that $\| f_m-g_m
\|_{C^0(U_{\alpha +\epsilon })}<m^{-1}$ and inevitably $\lim_m \|
\partial f_m(\xi )/
\partial \xi -\partial g_m(\xi )/\partial \xi \| =0$.
\par Let $V_{\xi }$ and $W_{\xi }$ be as above, then there exists
$\delta >0$ such that $(u(z)+{\tilde u}(z))/2<-\delta $ for each
$z\in U_{\alpha +\delta }\cap (V_{\xi }\setminus W_{\xi })$.
Therefore, there exists a branch of the ${\cal A}_p$ logarithm $Ln
(u)\in {\cal H}(U_{\alpha +\delta }\cap (V_{\xi }\setminus cl
(W_{\xi }))$ (see \S \S 3.7, 3.8 \cite{luoyst,luoystoc}). From
Theorems $4.2, 4.4$ it follows that each ${\cal A}_p$ holomorphic
Cousin problem over $U_{\alpha +\delta }$ has a solution. Hence $Ln
(u)=w_1-w_2$ for suitable $w_1\in {\cal H}(V_{\xi } \cap U_{\alpha
+\delta })$ and $w_2\in {\cal H}(U_{\alpha +\delta } \setminus cl
(W_{\xi })).$ Put $f:=u \exp (-w_1)$ in $U_{\alpha +\delta }\cap
V_{\xi }$ and $f:=\exp (-w_2)$ in $U_{\alpha +\delta } \setminus cl
(W_{\xi })$. Then $f\in {\cal H}(U_{\alpha +\delta })$ and $f(\xi
)=0$. In view of Inequality $(3)$ $f(z)\ne 0$ for each $\xi \ne z\in
cl(U_{\alpha })$. Verify now that $cl(U_{\alpha })$ is ${\cal
H}(M)$-convex. Consider $\xi \in M\setminus cl(U_{\alpha })$. Due to
\S 4.4 there exists a strictly plurisubharmonic exhausting
$C^{\omega }_{z,{\tilde z}}$-function $\psi $ for $M$ such that
$Crit (\psi )$ is discrete and $U_{\alpha }\subset G_{\psi (\xi )}$,
where $G_{\beta }:= \{ z\in M:$ $\psi (z)<\beta \} $ for $\beta \in
\bf R$. Considering shifts $\psi \mapsto \psi +const$ assume $d\psi
(z)\ne 0$ for each $z\in \partial G_{\psi (\xi )}$. From the proof
above it follows, that there exists $f\in {\cal H}(M)$ such that
$f(\xi )= 1$ and $|f(z)|<1$ for each $z\in cl(U_{\alpha })$.
\par {\bf 4.8.1. Remark.} With the help of Theorem $4.8$
it is possible to spread certain modifications of Theorems $3.2$ and
$3.5$ on ${\cal A}_p$ Stein manifolds.
\par {\bf 4.9. Theorem.} {\it Let $N$ be an ${\cal A}_s$
manifold with $1\le s\le \infty $ (where ${\cal A}_1:=\bf C$), then
for each $p$ with $s< p$ or $s\subset {\bf N}\subset p=\Lambda $,
$s\ne p$, there exists an ${\cal A}_p$ manifold $M$ and an ${\cal
A}_s$ holomorphic embedding $\theta : N\hookrightarrow M$.}
\par {\bf Proof.} Suppose $At (N)=\{ (V_a,\psi _a): a \in \Upsilon \} $
is any ${\cal A}_s$ holomorphic atlas of $N$, where $V_a$ is open in
$N$, $\bigcup_aV_a=N$, $\psi _a: V_a\to \psi _a(V_a)\subset {\cal
A}_s^n$ is a homeomorphism for each $a$, $n=dim_{{\cal A}_s}M\in \bf
N$, $\{ V_a: a \in \Upsilon \} $ is a locally finite covering of
$N$, $\psi _b\circ \psi _a^{-1}$ is a holomorphic function on $\psi
_a(V_a\cap V_b)$ for each $a, b\in \Lambda $ such that $V_a\cap
V_b\ne \emptyset $. Since ${\cal A}_p^n$ is normed, then it is
paracompact together with $M$ by Theorem 5.1.3 \cite{eng}. For each
${\cal A}_s$ holomorphic function $f$ on an open subset $V$ in
${\cal A}_s^n$ there exists an ${\cal A}_p$ holomorphic function $F$
on an open subset $U$ in ${\cal A}_p^n$ such that $\pi (U)=V$ and
$F|_V=f|_V$, where $\pi : {\cal A}_p^n\to {\cal A}_s^n$ is the
natural projection (see Proposition $3.13$ \cite{luoyst} and
analogously in the general case using local analyticity and a
locally finite covering of $V$).
\par Therefore, for each two charts $(V_a,\psi _a)$ and
$(V_b,\psi _b)$ with $V_{a,b}:= V_a\cap V_b\ne \emptyset $ there
exists $U_{a,b}$ open in ${\cal A}_p^n$ and an ${\cal A}_p$
holomorphic function $\Psi _{b,a}$ such that $\Psi _{b,a}|_{\psi
_a(V_{a,b})}=\psi _{b,a}|_{\psi _a(V_{a,b})}$, where $\psi
_{b,a}:=\psi _b\circ \psi _a^{-1}$, $\pi (U_{a,b})=\psi
_a(V_{a,b})$. Consider $Q:=\bigoplus_aQ_a$, where $Q_a$ is open in
${\cal A}_p^n$, $\pi (Q_a)=\psi _a(V_a)$ for each $a\in \Upsilon $.
The equivalence relation $\cal C$ in the topological space
$\bigoplus_a\psi _a(V_a)$ generated by functions $\psi _{b,a}$ has
an extension to the equivalence relation $\cal H$ in $Q$. Then
$M:=Q/\cal H$ is the desired ${\cal A}_p$ manifold with $At (M)=\{
(\Psi _a,U_a): a\in \Upsilon \} $ such that $\Psi _b\circ \Psi
_a^{-1}=\Psi _{b,a}$ for each $U_a\cap U_b\ne \emptyset $, $\Psi
_a^{-1}|_{\psi _a(V_a)}= \psi _a^{-1}|_{\psi _a(V_a)}$ for each $a$,
$\Psi _a^{-1}: Q_a\to U_a$ is the ${\cal A}_p$ homeomorphism.
Moreover, each homeomorphism $\psi _a: V_a\to \psi _a(V_a)\subset
{\cal A}_s^n$ has the ${\cal A}_p$ extension up to the homeomorphism
$\Psi _a: U_a\to \Psi _a(U_a)\subset {\cal A}_p^n$. The family of
embeddings $\eta _a: \psi _a(V_a)\hookrightarrow Q_a$ such that $\pi
\circ \eta _a=id$ together with $At(M)$ induces the ${\cal A}_s$
holomorphic embedding $\theta : N\hookrightarrow M$.
\par {\bf 4.10. Definition.} Let $M$ be an ${\cal A}_p$ manifold,
$2\le p\in \bf N$ or $p=\Lambda $. Suppose that for each chart
$(U_a,\phi _a)$ of $At(M)$ there exists an ${\cal A}_p$
superdifferentiable mapping $\Gamma : u\in \phi _a(U_a)\mapsto
\Gamma (u)\in L_q(X,X,X^*_q;{\cal A}_p)= L_q(X,X;X)$, where
$L_q(X^n,(X^*_q)^m;Y)$ denotes the space of all quasi-linear
mappings from $X^n\times (X^*_q)^m$ into $Y$ (that is, additive and
$\bf R$-homogeneous by each argument $x$ in $X$ or in $X^*_q$),
where $X$ and $Y$ are Banach spaces over ${\cal A}_p$, $X^*_q$
denotes the space of all additive $\bf R$-homogeneous functionals on
$X$ with values in ${\cal A}_p$ (see \S 4.1), $X^*_q=L_q(X;{\cal
A}_p)$. If $U_a\cap U_b\ne \emptyset $, let
\par $(1)\quad D(\phi _b\circ \phi _a^{-1}).\Gamma (\phi _a)=
D^2(\phi _b\circ \phi _a^{-1})+\Gamma (\phi _b)\circ (D(\phi _b\circ
\phi _a^{-1})\times D(\phi _b\circ \phi _a^{-1})).$ These $\Gamma
(\phi _a)$ are called the Christoffel symbols. Let ${\cal B}={\cal
B}(M)$ be a family of all ${\cal A}_p$ holomorphic vector fields on
$M$. For $M$ supplied with $\{ \Gamma (\phi _a): a \} $ define a
covariant derivation $(X,Y)\in {\cal B}^2\mapsto \nabla _XY\in \cal
B$:
\par $(2) \nabla _XY(u)=DY(u).X(u)+\Gamma (u)(X(u),Y(u))$,
where $X(u)$ and $Y(u)$ are the principal parts of $X$ and $Y$ on
$(U_a,\phi _a)$, $u=\phi _a (z)$, $z\in U_a$. In this case it is
said that $M$ possesses a covariant derivation.
\par {\bf 4.11. Remark.} Certainly for an ${\cal A}_p$ manifold
there exists a neighbourhood $V$ of $M$ in $TM$ such that $\exp :
V\to M$ is ${\cal A}_p$ holomorphic (see the real case in
\cite{kling}).
\par {\bf 4.12. Theorem.} {\it Let $f$ be an ${\cal A}_p$ holomorphic function
such that ${\hat f}$ is ${\cal A}_p$ (right) superlinear on a
compact ${\cal A}_p$ manifold $M$, where $2\le p\in \bf N$. Then $f$
is constant on $M$.}
\par {\bf Proof.} By the supposition of this theorem
$(f\circ \phi _b^{-1}){\hat .}$ is ${\cal A}_p$ (right) superlinear
for each chart $(U_b,\phi _b)$ of $M$. Since $M$ is compact and
$|f(z)|$ is continuous, then there exists a point $q\in M$ at which
$|f(z)|$ attains its maximum. Let $q\in U_b$, then $f\circ \phi
_b^{-1}$ is the ${\cal A}_p$ holomorphic function on $V_b:=\phi
_b(U_b)\subset {\cal A}_p^n$, where $dim_{{\cal A}_p}M=n$. Consider
a polydisk $V$ in ${\cal A}_p^n$ with the centre $y=\phi _b(q)$ such
that $V\subset V_b$. Put $g(w)=f\circ \phi _b^{-1} (y+(z-y)w)$,
where $w$ is the ${\cal A}_p$ variable. Then for each $z\in V$ there
exists $\epsilon _z>0$ such that the function $g(w)$ is ${\cal A}_p$
holomorphic on the set $W_z:=\{ w: w\in {\cal A}_p, |w|<1+\epsilon
_z \} $ and $|g(w)|$ attains its maximum at $w=0$. In view of
Theorem $3.15$ and Remark $3.16$ \cite{luoyst,luoystoc} $g$ is
constant on $W_z$, hence $f$ is constant on $U_b$. By the ${\cal
A}_p$ holomorphic continuation $f$ is constant on $M$.

\thanks{
Address: Sergey V. Ludkovsky, Mathematical Department, Brussels
University, Pleinlaan 1, Brussels, Belgium.\\
{\underline {Acknowledgment}}. The author thanks the Flemish Science
Foundation for support through the Noncommutative Geometry from
Algebra to Physics project and Professors Stefaan Caenepeel and Fred
van Oystaeyen for hospitality.}
\end{document}